\documentclass[titlepage,final,11pt]{article}
\usepackage{theorem}
\usepackage{enumerate}
\usepackage{array}
\usepackage[reqno]{amsmath}
\usepackage{amssymb}
\usepackage{mathtools}
\usepackage{latexsym}
\usepackage{makeidx}
\usepackage{fancybox}
\input epsf.sty
\usepackage[dvips]{graphicx,color}
\usepackage{showlabels}
\usepackage{subcaption}
\usepackage[dvipsnames]{xcolor}
\usepackage[normalem]{ulem}

\theoremstyle{change}
\newtheorem{proclaim}{PROCLAIM}[section]
\newtheorem{theorem}[proclaim]{Theorem}
\newtheorem{definition}[proclaim]{Definition}
\newtheorem{lemma}[proclaim]{Lemma}
\newtheorem{proposition}[proclaim]{Proposition}

\newtheorem{example}[proclaim]{Example}

\newtheorem{assumption}[proclaim]{Assumption}

\numberwithin{equation}{section}

\outer\def\proclaim #1. #2\par{\medbreak \noindent{\bf#1.\enspace}{\sl#2}\par
  \ifdim\lastskip<\medskipamount
  \removelastskip\penalty55\medskip\fi}
\def\state #1. { \noindent{\bf#1.\enspace}}
\def\algo #1. { \noindent{\bf#1.\enspace}}

\DeclareMathOperator{\dist}{dist}

\DeclareMathOperator{\dom}{dom}

\newcommand{\comp}{\,{\raise 1pt \hbox{$\scriptstyle\circ$}}\,}

\newcommand{\reals}{\mathbb{R}}

\newcommand{\natnums}{{{\rm l} \kern -.13em {\rm N} }}
\newcommand{\nats}{\mathbb{N}}
\newcommand{\snats}{{I\kern -.29em N}}
\newcommand{\rats}{{Q\kern -.64em \raise 1pt \hbox{$\scriptstyle |$}\;\,}}
\newcommand{\srats}
	{{Q\kern -.56em \raise 1.2pt \hbox{$\scriptscriptstyle /$}\,}}
\newcommand{\ints}{Z\kern -.46em Z}
\newcommand{\ball}{\mathbb{B}}
\newcommand{\pluss}{\hskip1pt \raise1pt\vbox{\hrule width6pt \vskip1pt \hrule
                    width6pt} \kern-4pt{\lower1pt\hbox{\vrule height6pt
		    \kern1pt\vrule height6pt}}\hskip5pt}
\newcommand{\eop}
	{\hfill{$\vcenter{\hrule height1pt \hbox{\vrule width1pt height5pt
   	 \kern5pt \vrule width1pt} \hrule height1pt}$} \medskip}

\newcommand{\setd}{{ d \kern -.15em l}}
\newcommand{\hatsetd}{ d \hat{\kern -.15em l }}

\renewcommand{\epsilon}{\varepsilon}
\renewcommand{\phi}{\varphi}

\newcommand{\tto}{\;{\lower 1pt \hbox{$\rightarrow$}}\kern -12pt
           \hbox{\raise 2.5pt \hbox{$\rightarrow$}}\;}
\newcommand{\overto}[1]{\,{\raise 0pt\hbox{$\rightarrow$}}\kern -9pt
     \hbox{\lower 3pt \hbox{$\scriptscriptstyle#1$}}\hskip6pt}
\newcommand{\underto}[1]{\,{\lower 1pt\hbox{$\rightarrow$}}\kern -9pt
     \hbox{\raise 4pt \hbox{$\,\scriptscriptstyle#1$}}\hskip7pt}
\newcommand{\bigoverto}[1]{{\raise 0pt\hbox{$\,\longrightarrow$}}\kern -16pt
     \hbox{\lower 3pt \hbox{$\scriptscriptstyle#1$}}\hskip4pt}
\newcommand{\bigunderto}[1]{\,{\lower 1pt\hbox{$\longrightarrow$}}\kern -16pt
     \hbox{\raise 4pt \hbox{$\,\scriptscriptstyle#1$}}\hskip6pt}
\newcommand{\bigbigto}[2]{\,{\raise 0pt\hbox{$\,\longrightarrow$}}\kern -16pt
     \hbox{\lower 3pt \hbox{$\scriptscriptstyle#2$}}\kern -10pt
     \hbox{\raise 4pt \hbox{$\,\scriptscriptstyle#1$}}\hskip7pt}
\newcommand{\downto}{{\raise 1pt \hbox{$\scriptscriptstyle \,\searrow\,$}}}
\newcommand{\upto}{{\raise 1pt \hbox{$\scriptscriptstyle \,\nearrow\,$}}}

\newcommand{\notimply}
	{\quad\hbox{$\Longrightarrow \kern -14pt {/}$}\hskip6pt\quad}

\newcommand{\lto}{\,{\lower 1pt\hbox{$\rightarrow$}}\kern -10pt
     \hbox{\raise 4pt \hbox{$\, \scriptstyle l$}}\hskip7pt}
\newcommand{\eto}{\,{\lower 1pt\hbox{$\rightarrow$}}\kern -10pt
     \hbox{\raise 4pt \hbox{$\, \scriptstyle e$}}\hskip7pt}
\newcommand{\hto}{\,{\lower 1pt\hbox{$\rightarrow$}}\kern -11pt
     \hbox{\raise 4pt \hbox{$\, \scriptstyle h$}}\hskip7pt}
\newcommand{\pto}{\,{\lower 1pt\hbox{$\rightarrow$}}\kern -11pt
     \hbox{\raise 4.5pt \hbox{$\, \scriptstyle p$}}\hskip7pt}
\newcommand{\cto}{\,{\lower 1pt\hbox{$\rightarrow$}}\kern -11pt
     \hbox{\raise 4pt \hbox{$\, \scriptstyle c$}}\hskip7pt}
\newcommand{\gto}{\,{\lower 1pt\hbox{$\rightarrow$}}\kern -11pt
     \hbox{\raise 4.5pt \hbox{$\, \scriptstyle g$}}\hskip7pt}
\newcommand{\sto}{\,{\lower 1pt\hbox{$\rightarrow$}}\kern -10pt
     \hbox{\raise 4pt \hbox{$\, \scriptstyle s$}}\hskip7pt}
\newcommand{\awto}{\,{\lower 1pt\hbox{$\rightarrow$}}\kern -15pt
     \hbox{\raise 4pt \hbox{$\, \scriptstyle aw$}}\hskip7pt}
\def\Nto{\,{\raise 1pt\hbox{$\rightarrow$}}\kern -13pt
     \hbox{\lower 3pt \hbox{$\, \scriptstyle N$}}\hskip7pt}
\def\Cto{\,{\raise 1pt\hbox{$\rightarrow$}}\kern -14pt
     \hbox{\lower 3pt \hbox{$\, \scriptstyle C$}}\hskip7pt}
\def\fto{\,{\raise 1pt\hbox{$\rightarrow$}}\kern -14pt
     \hbox{\lower 3pt \hbox{$\, \scriptstyle f$}}\hskip7pt}

\newcommand{\low}[1]{{\lower1pt \hbox{$\scriptstyle #1$}}}
\newcommand{\loww}[1]{{\lower2pt \hbox{$\scriptstyle #1$}}}
\newcommand{\high}[1]{{\raise1pt \hbox{$\scriptstyle #1$}}}

\newcommand{\nlim}{\mathop{\rm lim}\nolimits}
\newcommand{\nliminf}{\mathop{\rm liminf}\nolimits}
\renewcommand{\liminf}{\mathop{\rm liminf}}

\newcommand{\nlimsup}{\mathop{\rm limsup}\nolimits}
\newcommand{\ninf}{\mathop{\rm inf}\nolimits}

\newcommand{\nnmin}{\mathop{\rm minimize}}

\newcommand{\nargmin}{\mathop{\rm argmin}\nolimits}

\newcommand{\lwdy}[2]{\mathrel{\mathop
        {\raisebox{0.1ex}{\null$#1$}}{\hbox{\kern -1.0em
	{\raisebox{-0.8ex}{$\scriptstyle{\;\to #2}$}}}}}}
\newcommand{\lwwdy}[2]{\mathrel{\mathop
        {\raisebox{0.2ex}{\null$#1$}}{\hbox{\kern -1.0em
	{\raisebox{-1.1ex}{$\scriptstyle{\;\to #2}$}}}}}}
\newcommand{\slwdy}[2]{\scriptsize{{\mathrel{\mathop
        {\raisebox{0.1ex}{\null$#1$}}{\hbox{\kern -1.0em
	{\raisebox{-0.8ex}{$\scriptstyle{\;\to #2}$}}}}}}}}
\newcommand{\slwwdy}[2]{\scriptsize{{\mathrel{\mathop
        {\raisebox{0.2ex}{\null$#1$}}{\hbox{\kern -1.0em
	{\raisebox{-1.1ex}{$\scriptstyle{\;\to #2}$}}}}}}}}

\definecolor{lightgray}{gray}{0.75}
\definecolor{myred}{rgb}{0.55,0,0}
\definecolor{myblue}{rgb}{0,0,0.5} % hex: #00007f
\definecolor{mygreen}{rgb}{0,0.5,0} % hex: #00007f
\definecolor{purple}{rgb}{0.5,0,0.5} % hex: #00007f
\definecolor{turq}{rgb}{0,0.805,0.816} % hex: #00007f
\definecolor{maroon}{rgb}{0.51,0,0}
\definecolor{MAROON}{rgb}{0.51,0,0}
\definecolor{redor}{rgb}{0.78,0.078,0.078}
\definecolor{dgreen}{rgb}{0,0.3,0}
\newcommand{\bblue}{\color{black}}
\newcommand{\blue}{\color{black}}

\newcommand{\bcdot}{\,{\raise .2ex \hbox{$\centerdot$}}\,}

\begin{document}

%% FOOTNOTE MARK: *, \dag ... instead of  1, 2, ...
%\renewcommand{\thefootnote}{\fnsymbol{footnote}}

\begin{center}
\begin{large}
{\bf On Stability in Optimistic Bilevel Optimization}
\smallskip
\end{large}
\vglue 0.5truecm
\begin{tabular}{cc}
  \begin{large} {\sl Johannes O. Royset} \end{large}\\
  Daniel J. Epstein Department of Industrial and Systems Engineering\\
  University of Southern California
\end{tabular}

\vskip 0.2truecm

\end{center}

\vskip 0.6truecm

\noindent {\bf Abstract}. Solutions of bilevel optimization problems tend to suffer from instability under changes to problem data. In the optimistic setting, we construct a lifted formulation that exhibits desirable stability properties under mild assumptions that neither invoke convexity nor smoothness. The upper- and lower-level problems might involve integer restrictions and disjunctive constraints. In a range of results, we invoke at most pointwise and local calmness for the lower-level problem in a sense that holds broadly. The lifted formulation is computationally attractive with structural properties being brought out and an outer approximation algorithm becoming available.

\vskip 0.2truecm

%% KEDWORDS, AMS CLASSIFICATION, DATE
\halign{&\vtop{\parindent=0pt
   \hangindent2.5em\strut#\strut}\cr
{\bf Keywords}: Bilevel optimization, optimistic formulation, stability, epi-convergence.
   %% second line of key words if necessary
                         \cr
%{\bf AMS Classification}: \quad   %% AMS numbers xxXxx, xxXxx
%                          \cr\cr

{\bf Date}:\quad \ \today \cr}

\baselineskip=15pt

\section{Introduction}\label{sec:intro}

Bilevel problems tend to be ill-posed in the sense that small changes to problem data may cause disproportional large changes in solutions. In fact, ill-posedness can be expected as it is well-known from variational analysis that the set of minimizers to a function is stable in a meaningful sense only under strong assumptions. Even when the lower-level problem is seemingly ``nice,'' a small change in a lower-level constraint may have a large effect on the solution of a bilevel problem; \cite{BeckBienstockSchmidtThurauf.23} furnishes an examples with a linear lower-level objective function, a convex lower-level feasible set satisfying the Slater condition, and a unique lower-level solution satisfying a strict complementarity condition for all upper-level feasible decisions. Additional examples appear in \cite{CarusoLignolaMorgan.20}. Since problem data is almost always unsettled in practice, it is clear that blindly solving a bilevel problem without any safeguards may not produce meaningful decisions. In this paper, we provide a lifted formulation of bilevel problems that is fundamentally well-posed in the sense that discrepancies in the problem data still lead to useful decisions and bounds in the limit as the discrepancies vanish. While the lifted formulation may involve a large number of inequality constraints, it preserves structural properties of the upper- and lower-level problems and thus can be approaches by many well-known algorithmic strategies including those from semi-infinite optimization.  

We consider the following bilevel problem. For a nonempty set $X\subset\reals^n$ and a function $f:\reals^n\times\reals^m \to (-\infty,\infty]$, the upper-level problem is to minimize $f(x,y)$ over $x\in X$, where the lower-level problem furnishes $y\in \reals^m$. For nonempty sets $Y\subset\reals^m$ and $D\subset\reals^q$, a function $g:\reals^n\times\reals^m\to \reals$, and a mapping $H:\reals^n\times\reals^m\to \reals^q$, the lower-level problem amounts to minimizing $g(x,y)$ over $y\in Y$ while satisfying $H(x,y) \in D$. Since $D$ could be the product set of $\{0\}$ and $(-\infty, 0]$, the inclusion $H(x,y) \in D$ may represent any finite number of equality and inequality constraints. Disjunctive constraints emerge from nonconvex $D$. Throughout, we adopt the {\em optimistic} perspective. (We refer to \cite{LiuFanChenZheng.18,LiuFanChenZheng.20} for the pessimistic perspective and to \cite{AboussororLoridan.95,CarusoLignolaMorgan.20} for other alternatives.) This leads to the {\em bilevel problem}   
\begin{align}
\mbox{(P)}~~~~~~~  &\nnmin_{x\in X, y\in \reals^m}  ~ f(x,y)\nonumber\\
& \mbox{ subject to } ~~ y \in \tau\mbox{-}\nargmin_{z\in Y} \big\{ g(x,z) ~\big|~ H(x,z) \in D\big\},\label{eqn:argminconstr}
\end{align}
where $\tau \in [0,\infty)$ allows us to also consider near-optimal solutions of the lower-level problem by sometimes setting $\tau > 0$; recall that $\epsilon\mbox{-}\nargmin_{z\in Z} \psi(z) = \{\bar z \in Z ~|~ \psi(\bar z) \leq \inf_{z\in Z} \psi(z) + \epsilon, \psi(\bar z) < \infty\}$ for any extended real-valued function $\psi$, set $Z$, and tolerance $\epsilon \in [0,\infty)$. Since $f$ is extended real-valued, (P) may have {\em implicit constraints} in addition to the {\em explicit constraints} $x\in X$ and \eqref{eqn:argminconstr}. Thus, its {\em feasible set} $\{(x,y)\in X\times \reals^m \, |\, f(x,y)<\infty, \mbox{\eqref{eqn:argminconstr} holds}\}$ may impose coupling constraints between $x$ and $y$. 

To model unsettled data in (P), we consider $f^\nu:\reals^n\times\reals^m\to (-\infty, \infty]$, $g^\nu:\reals^n\times\reals^m\to \reals$, and $H^\nu:\reals^n\times\reals^m \to \reals^q$ as approximations of $f, g$, and $H$, respectively. We examine the situation as these approximations become more accurate as indicated by higher values of the superscript $\nu\in \nats = \{1, 2, \dots\}$. As alluded to above, replacing $f,g,H$ in (P) by $f^\nu,g^\nu,H^\nu$ produces the {\em naive problem} 
\begin{align*}
\mbox{(NP)}^\nu~~~~~~~  &\nnmin_{x\in X, y\in \reals^m}  ~ f^\nu(x,y)\\
& \mbox{ subject to } ~~ y \in \tau\mbox{-}\nargmin_{z\in Y} \big\{ g^\nu(x,z) ~\big|~ H^\nu(x,z) \in D\big\},
\end{align*}
which may have solutions far from those of (P) even if the corresponding functions are close, for example in the sense of the sup-norm on $X\times Y$; Section \ref{sec:example} furnishes concrete instances. In this paper, we propose the following {\em lifted problem} as an alternative: 
\begin{align}
\mbox{(P)}^\nu~~~~~~~~~~~~~  \nnmin_{x, y, u, \alpha, \lambda} ~ f^\nu(x,y) & + \sigma^\nu \|u\| - \theta^\nu \alpha\nonumber\\
  \mbox{subject to } H^\nu(x,y) + u & \in D\nonumber\\
  g^\nu(x,y) + \alpha & \leq g^\nu(x,z) + \lambda \dist\big(H^\nu(x,z), D\big) + \tau^\nu~~\forall z\in Y^\nu \label{eqn:PnuConstr_mu}\\
  x \in X, ~y & \in Y, ~ u \in \reals^{q}, ~ \alpha \in (-\infty, 0], ~ \lambda \in [0,  \bar \lambda^\nu].\nonumber 
\end{align}
The problem involves the nonnegative parameters $\sigma^\nu, \theta^\nu, \bar\lambda^\nu$, and $\tau^\nu$. The two first ones are penalty parameters for the auxiliary decision variables $u\in \reals^q$ and $\alpha\in \reals$. (Squared and other penalty terms are also possible, but we omit the details.) The auxiliary decision variable $\lambda$ is bounded from above by $\bar\lambda^\nu$. The last parameter, $\tau^\nu$, is a counterpart to the tolerance $\tau$ for the lower-level problem in (P). Given any point $a\in \reals^q$, $\dist(a, D) = \inf_{a' \in D} \|a - a'\|$ is the {\em point-to-set distance} under an {\em arbitrary norm} $\|\cdot\|$ on $\reals^q$. The norm underpinning $\dist(\cdot,D)$ in (P)$^\nu$ matches the norm used in the term $\sigma^\nu\|u\|$ of the objective function. For example, 
if $D = (-\infty, 0]^q$, then $\dist(H^\nu(x,z), D)=\sum_{j = 1}^q \max\{0, h^\nu_j(x,z)\}$ when $\|\cdot\| = \|\cdot\|_1$ and $h_j^\nu$ is the $j$th component of $H^\nu$.

Lastly, (P)$^\nu$ involves a set $Y^\nu \subset Y$, which could coincide with $Y$. Since $Y$ tends to have a high or infinite cardinality, the possibility of considering a smaller set $Y^\nu$ in \eqref{eqn:PnuConstr_mu} is practically beneficial. 

At first glance, (P)$^\nu$ appears only loosely related to (P). However, the next section establishes that solving (P)$^\nu$,  despite the use of approximating data $(f^\nu, g^\nu, H^\nu, Y^\nu)$, bounds the minimum value of (P) in the limit as $\nu\to \infty$ under mild assumptions. We also show that (P)$^\nu$ can produce solutions that in the limit are feasible and even optimal for (P). This takes place in settings where (NP)$^\nu$ fails to produce such convergence properties.  

In addition to have desirable limiting properties as $\nu\to\infty$, (P)$^\nu$ is constructed with computations in mind. While there are several nontrivial challenges associated with solving (P)$^\nu$, it is apparent that any structural properties in $g^\nu$ and $H^\nu$ such as convexity or linearity are being brought forward. Moreover, $Y^\nu$ does not depend on $x$, which facilitates the use of outer approximations as seen in Section \ref{sec:computational}. 

There is an extensive literature on bilevel problems. The monographs \cite{Dempe.02,DempeKalashnikovPerezvaldesKalashnykova.15} provide an overview of the main approaches and applications; see also the collection \cite{DempeZemkoho.20}, especially \cite{Dempe.20} for a recent list of references, and \cite{KleinertLabbeLjubicSchmidt.21} for details about computational methods in the setting of mixed-integer programming. Modern applications in machine learning elevates the importance of bilevel optimization even higher; see, e.g., \cite{OkunoTakeda.20,KwonKwonWrightNowak.23a,ChenXiaoBalasubramanian.24,KwonKwonWrightNowak.23b}. When the lower-level problem is changed from an optimization problem to an equilibrium problem, possibly triggered by optimality conditions of the lower-level problem, we obtain further extensions as discussed in \cite{LuoPangRalphWu.96,YeZhuZhu.97}; see also the monograph \cite{LuoPangRalph.83}. 

%Optimality conditions arising in this setting may be enforced using penalty terms as in \cite{MarcotteZhu.96}. 

%A summary of extensions to models accounting for uncertainty appears in \cite{BeckLjubicSchmidt.23}. 

Our lifted formulation (P)$^\nu$ is most closely related to approaches in which the argmin constraint \eqref{eqn:argminconstr} is reformulated by a constraint involving the min-value function for the lower-level problem, i.e., the {\em value-function approach} \cite{Outrata.90,ChenFlorian.91,YeZhu.95,KwonKwonWrightNowak.23b}, as well as related methods such as the two-level value-function approach \cite{DempeMordukhovichZemkoho.12,DempeMordukhovichZemkoho.19} and techniques bounding the min-value function \cite{MitsosLemonidisBarton.08,Mitsos.10,DominguezPistikopoulos.10,KoppeQuayranneRyan.10}. In our notation, the value-function approach leverages the trivial fact: 
For $x\in\reals^n$ and $\tau \in [0,\infty)$, one  has
\begin{align}
& y\in \tau\mbox{-}\nargmin_{z\in Y} \big\{g(x,z) \,\big|\, H(x,z) \in D\big\} ~~~~~~~\Longleftrightarrow~~\label{eqn:level-reform}\\
& y\in Y, ~~H(x,y) \in D, ~~g(x,y) \leq  \inf_{z\in Y} \big\{g(x,z) \,\big|\, H(x,z) \in D\big\} + \tau.\nonumber 
\end{align}  
Although just for equality and inequality constraints in the lower-level problem, \cite{YeZhu.95} utilizes this equivalence as it introduces {\em partial calmness} as a qualification that allows the constraint bounding $g(x,y)$ in \eqref{eqn:level-reform} to be moved into the upper-level objective function using exact penalization. That paper identifies uniformly weak sharp minima for the lower-level problem as a sufficient condition for partial calmness; see also the treatment in \cite[Chapter 6]{Mordukhovich.18}. In machine learning applications, \cite{KwonKwonWrightNowak.23b,LiuLiuZengZhangZhang.23} also leverage such penalization. {\blue The former paper omits the inclusion $H(x,z) \in D$, while the latter uses a standard barrier functions to address such constraints but these become problematic in the absence of strict lower-level feasibility.} Smoothing of the min-value function is also viable \cite{LinXuYe.14,AlcantaraTakeda.24}, at least when there is no $x$-dependence in the lower-level constraints.

While also utilizing \eqref{eqn:level-reform}, we proceed in an manner distinct from these efforts. We move the constraint $H(x,z) \in D$ into the {\em lower-level} objective function and sometimes, but not always, use calmness to ensure exact penalization. After also introducing approximating data, this step leads to \eqref{eqn:PnuConstr_mu}. {\blue (We prefer exact penalization over log-barrier type penalties of $H(x,z) \in D$ because of difficulties for the latter approach in the absence of strict feasibility.)} The approach is well-aligned with the classical connection between calmness and exact penalization; see \cite[Prop. 6.4.3]{Clarke.83} and \cite{Burke.91,Burke.91b}. Our use of calmness is distinct from partial calmness of \cite{YeZhu.95} by placing the focus on the lower-level $z$-vector and not the upper-level $x$-vector; see Subsection \ref{subsec:calmness} for details. {\blue The resulting lifted problem (P)$^\nu$ is computationally attractive and amenable to further algorithmic developments as we see below.} 

In contrast, the reformulation \eqref{eqn:level-reform} leads to the more challenging generalized semi-infinite constraints $g(x,y) \leq g(x,z) + \tau$ for each $z\in Y$ that also satisfies $H(x,z) \in D$. After introducing a feasibility tolerance, \cite{TsoukalasRustemPistikopoulos.09} approximates these constraints via a min-max-min problem. Further algorithmic efforts to handle generalized semi-infinite constraints include the use of generalized Benders' decomposition  
\cite{BolusaniConiglioRalphsTahernejad.20,BolusaniRalph.22}. The value-function approach also underpins \cite{LozanoSmith.17} and its algorithm for specially structured bilevel problems. One may attempt to ``track'' how the lower-level optimal solutions change as a function of the upper-level vector $x$; see, e.g., \cite{Dempe.00,SteinStill.02,FaiscaDuaRustemSaraivaPistikopoulos.07}. Here, regularization may mitigate some difficulties. Regularization comes in the form of a norm-squared penalty in the lower-level problem \cite{Dempe.00} and by allowing near-optimal solutions in the lower-level problem \cite{LinXuYe.14}. The multitude of lower-level optimal solutions may also be addressed via multi-objective optimization techniques \cite{Zemkoho.16}.

Stability results for bilevel optimization problems extend at least back to \cite{LoridanMorgan.84,LoridanMorgan.85,LoridanMorgan.89}. Most significantly, \cite{LignolaMorgan.97} (see also \cite{CarusoLignolaMorgan.20}) leverages notions of variational convergence in the study of bilevel problems under perturbations and in particular convergence of the min-value function of the lower-level problem. It is recognized that a regularization of the lower-level problem by allowing near-optimal solutions improves stability properties. Still, \cite{LignolaMorgan.97} invokes arguably strong assumptions including the existence of interior points for the lower-level feasible set and convexity and quasiconvexity for certain components; examples illustrate how these assumptions appear nearly unavoidable. We continue in this tradition, but define (P)$^\nu$ as a way around the challenges identified in \cite{LignolaMorgan.97} and elsewhere. The added variables in (P)$^\nu$ and the penalization of lower-level constraints to produce \eqref{eqn:PnuConstr_mu} appear to be the key to avoid strong assumptions. {\blue A direct use of the reformulation \eqref{eqn:level-reform} with approximating data would require additional assumptions due to the instability of minimum values under perturbations.} 

There are several issues that are beyond the scope of the paper. We solely focus on the optimistic formulation (P) and thus omit a treatment of its many alternatives. We concentrate on recovering global solutions of (P) by computing near-optimal solutions of (P)$^\nu$. Local and stationary solutions are more challenging to treat for several reasons including some listed in \cite{DempeMordukhovichZemkoho.12}. We shy away from quantifying the rate of convergence, which is difficult for methods relying on constraint penalization. 

After a summary of notation and terminology, we continue in Section \ref{sec:main} with the main theorems. Section \ref{sec:example} gives examples. Section \ref{sec:computational} discusses computational strategies for solving (P)$^\nu$ and establishes the validity of an outer approximation algorithm. Section \ref{sec:proofs} provides supplementary results and proofs.

\smallskip

\noindent{\bf Notation and Terminology}. The {\em Euclidean  ball} is $\ball(x,\rho) = \{x'\in \reals^n~|~\|x-x'\|_2\leq \rho\}$. For any set $C$, $\iota_C(x) = 0$ if $x\in C$ and $\iota_C(x) = \infty$ otherwise. A sequence of sets $C^\nu\subset\reals^n$ converging to $C\subset\reals^n$ in the sense of {\em Painlev\'{e}-Kuratowski} takes place when $C$ is closed and $\dist(x,C^\nu)\to \dist(x,C)$ for all $x\in\reals^n$, and this is denoted by $C^\nu\sto C$. For any function $\psi:\reals^n\to [-\infty,\infty]$, $\dom \psi = \{x \in \reals^n\,|\,\psi(x)<\infty\}$ is its {\em domain}, $\inf \psi = \inf_{x\in \reals^n} \psi(x)$ is its {\em minimum value}, and $\epsilon$-$\nargmin \psi = \{x\in \dom \psi~|~\psi(x) \leq \inf \psi + \epsilon\}$ is its set of $\epsilon$-{\em minimizers}. It is {\em lower semicontinuous} (lsc) if $\liminf \psi(x^\nu)\geq \psi(x)$ whenever $x^\nu\to x$ and it is {\em continuous relative to a set} $C$ if $\psi(x^\nu) \to \psi(x)$ whenever $x^\nu\in C\to x\in C$. The functions $\psi^\nu:\reals^n\to [-\infty, \infty]$ {\em epi-converge} to $\psi$, written $\psi^\nu\eto \psi$, when
\begin{align}
  &\forall x^\nu \to x, ~~\nliminf \psi^\nu(x^\nu) \geq \psi(x)\label{eqn:liminf}\\
  &\forall x, ~\exists x^\nu\to x \mbox{ with } \nlimsup \psi^\nu(x^\nu)\leq \psi(x).\label{eqn:limsup}
\end{align}
We say that $\{\psi^\nu, \nu\in\nats\}$ is {\em tight} if for all $\epsilon>0$, there exist compact $B \subset \reals^n$ and $\bar\nu\in \nats$ such that
$\ninf_{x\in B} \psi^\nu(x) \leq \inf \psi^\nu + \epsilon$ for all $\nu\geq \bar \nu$.

For any subsequence $N\subset \nats$, we write $x^\nu \Nto x$ to indicate convergence of $\{x^\nu, \nu\in N\}$ to $x$. The {\em minimum value} of a minimization problem is the infimum of its objective function values at points satisfying its explicit constraints. (It would be equal to infinity if there are no such points.) Specifically, 
\[
\mbox{$\mathfrak{m}$ is the {\em minimum value} of (P)~~ and ~~$\mathfrak{m}^\nu$ is the {\em minimum value} of (P)$^\nu$}. 
\]
An $\epsilon$-{\em optimal solution} of (P) is a pair $(x,y)\in X\times \reals^m$ satisfying \eqref{eqn:argminconstr}, $f(x,y) < \infty$, and $f(x,y) \leq \mathfrak{m} + \epsilon$. A $0$-optimal solution of (P) is called an {\em optimal solution} of (P). Analogous definitions apply for (P)$^\nu$.

\section{Main Results}\label{sec:main}

We develop in this section a series of results about the relationship between (P)$^\nu$ and (P) under gradually stronger assumptions starting with the following basic one. 

\begin{assumption}{\rm (basic assumption).}\label{ass:basic} The following hold:
\begin{enumerate}[(a)]

\item The sets $X\subset\reals^n$, $Y\subset\reals^m$, and $D\subset\reals^q$ are nonempty and closed. 

\item The sets $Y^\nu$, $\nu\in\nats$, are nonempty subsets of $Y$ and $Y^\nu\sto Y$. 

\item For a vanishing sequence $\{\delta^\nu, \nu\in\nats\}$, the functions $g,g^\nu:\reals^n\times\reals^m\to \reals$ satisfy 
\begin{equation*}\label{eqn:errorratesg}
\sup_{x\in X, y\in Y} \big| g^\nu(x,y) - g(x,y) \big| \leq \delta^\nu  ~~\forall \nu\in\nats
\end{equation*}
and $g$ is continuous relative to $X\times Y$. 

\item For a vanishing sequence $\{\eta^\nu, \nu\in\nats\}$, the mappings $H,H^\nu:\reals^n\times\reals^m\to \reals^q$ satisfy
\begin{equation}\label{eqn:errorratesH}
\sup_{x\in X, y\in Y} \Big| \dist\big( H^\nu(x,y), D\big) - \dist\big( H(x,y), D\big) \Big| \leq \eta^\nu ~~\forall \nu\in\nats
\end{equation}
and $H$ is continuous relative to $X\times Y$. 

\item For $x\in X$ and $y\in Y$, the functions $f,f^\nu:\reals^n\times\reals^m\to (-\infty,\infty]$ satisfy 
\begin{align*}
\nlimsup f^\nu(x,y) & \leq f(x,y)\\
\nliminf f^\nu(x^\nu,y^\nu) & \geq f(x,y) ~\mbox{ whenever } ~x^\nu\in X\to x, ~y^\nu\in Y\to y. 
\end{align*}

\item The nonnegative parameters $\sigma^\nu, \theta^\nu, \bar\lambda^\nu$ satisfy $\bar\lambda^\nu\to\infty$, $\sigma^\nu \eta^\nu \to 0$, $\theta^\nu\bar\lambda^\nu\eta^\nu\to 0$, and $\theta^\nu \delta^\nu\to 0$.

\end{enumerate}

\end{assumption}

Assumption \ref{ass:basic} is mild and captures a vast array of bilevel problems and their perturbations. Both the upper- and lower-level problem may involve integer restrictions and the lower-level problem may involve disjunctions as $D$ could be nonconvex. There are no requirements related to convexity and smoothness of $f,f^\nu,g,g^\nu,H,H^\nu$ or ``well-postedness'' of the lower-level problem such as the Slater condition or strict complementarity.  

Assumption \ref{ass:basic}(f) prescribes that the parameter $\bar\lambda^\nu$ should not grow too fast as compared to the vanishing approximation error in \eqref{eqn:errorratesH}. While not required by Assumption \ref{ass:basic}(f), we sometimes below let $\sigma^\nu,\theta^\nu\to \infty$. Then, for example, when $\delta^\nu = \eta^\nu = \nu^{-\beta}$, with $\beta \in (0,\infty)$, one can set $\theta^\nu = \bar\lambda^\nu = \nu^{\beta/3}$ and $\sigma^\nu = \nu^{\beta/2}$ and still satisfy Assumption \ref{ass:basic}(f). Regardless, the exact values of $\delta^\nu$ and $\eta^\nu$ can remain unknown as Assumption \ref{ass:basic}(f) only relies on their rates of decay.   

Proofs omitted in this section are given in Section \ref{sec:proofs}.

\subsection{Bounds and Feasible Solutions}\label{subsec:withoutcalmness}

We start by developing conditions under which (P)$^\nu$ leads to bounds on the minimum value of (P) as well as convergence to a feasible solution of (P). 

\begin{theorem}{\rm (bounds and feasible solutions).}\label{thm:withoutcalmness}
Under Assumption \ref{ass:basic}, the following hold: 
\begin{enumerate}[(a)]

\item If $\{\theta^\nu, \nu\in\nats\}$ is bounded away from zero, $\theta^\nu|\tau^\nu - \tau'|\to 0$ for some $\tau'>\tau$, and $Y$ is compact, then 
\[
\nlimsup \mathfrak{m}^\nu \leq \mathfrak{m}.
\]

\item Suppose there is a compact set $B\subset \reals^{n+m+q+2}$ such that for each $\nu\in\nats$ there exists an optimal solution of {\rm (P)}$^\nu$ in $B$. If  $\theta^\nu|\tau^\nu-\tau|\to 0$ and $\sigma^\nu,\theta^\nu\to \infty$, then  
    \[
    \nliminf \mathfrak{m}^\nu \geq \mathfrak{m},
    \]
    with $\{\mathfrak{m}^\nu, \nu\in\nats\}$ actually converging to a real number when $\nliminf \mathfrak{m}^\nu$ is finite. 

\item Suppose that $\theta^\nu|\tau^\nu-\tau|\to 0$, $\sigma^\nu$, $\theta^\nu\to \infty$, and there are $x\in X$, $y\in Y$, and $\lambda \in [0,\infty)$ such that 
    \begin{equation}\label{eqn:extraass}
f(x,y) < \infty,~ H(x,y)\in D, ~  g(x,y) \leq \inf_{z\in Y} g(x,z) + \lambda \dist\big(H(x,z), D\big) + \tau. 
    \end{equation}
    For every $\nu\in\nats$, let $(x^\nu, y^\nu, u^\nu, \alpha^\nu, \lambda^\nu)$ be an $\epsilon^\nu$-optimal solution of {\rm (P)}$^\nu$. If $\{\epsilon^\nu, \nu\in\nats\}$ is bounded and there is a subsequence $N\subset\nats$ and $(\hat x, \hat y, \hat u, \hat \alpha, \hat \lambda)$ such that 
    \[
    (x^\nu, y^\nu, u^\nu, \alpha^\nu, \lambda^\nu)\Nto (\hat x, \hat y, \hat u, \hat \alpha, \hat \lambda), 
    \]
    then $(\hat x, \hat y)$ is feasible in {\rm (P)}. 
 
 If $Y$ is compact, $\tau \in (0,\infty)$, and $f$ is real-valued, then the existence of $x\in X$, $y\in Y$, and $\lambda\in [0,\infty)$ satisfying \eqref{eqn:extraass} is guaranteed by the existence of a feasible point in {\rm (P)}.  
\end{enumerate}

\end{theorem}

Part (a) of the theorem asserts that the minimum value of (P)$^\nu$ furnishes in the limit a lower bound on the minimum value of (P), which in turn may help with assessing the quality of any incumbent solution of (P). Since $\nlimsup \mathfrak{m}^\nu$ is nondecreasing as $\tau^\nu$ decreases, it is beneficial to select $\tau'$ only slightly above $\tau$ (specified by (P) as a nonnegative number). The compactness of $Y$ stems from Lemma \ref{lemma:convmu} below, where it ensures that a penalization of the constraint $H(x,z) \in D$ is valid in some sense. After that lemma, we provide an example of what may go wrong in the absence of compactness. 

Theorem \ref{thm:withoutcalmness}(a) allows $\sigma^\nu$ and $\theta^\nu$ to remain bounded. Since $\mathfrak{m}^\nu$ is nondecreasing in these parameters, we are nevertheless motivate to choose higher values. An examination of the proof of Theorem \ref{thm:withoutcalmness}(a) shows that the liminf-condition in Assumption \ref{ass:basic}(e) is superfluous. However, approximating functions $f^\nu$ that remain well ``below'' $f$ would naturally produce a loose lower bound via $\nlimsup \mathfrak{m}^\nu$. 

An example illustrates that a lower bound on $\mathfrak{m}$, parallel to Theorem \ref{thm:withoutcalmness}(a), cannot be guaranteed by (NP)$^\nu$. 

\begin{example}{\rm (lower bound and role of $u$).}\label{ex:infeaslower}
For $g(x,y) = g^\nu(x,y) = 0$, $Y = Y^\nu = [-1,1]$, $D = (-\infty, 0]$, $H(x,y) = y^2$, and $H^\nu(x,y) = y^2 +1/\nu$, we obtain 
\begin{align*}
&\tau\mbox{-}\nargmin_{y\in Y}\big\{g(x,y)~\big|~H(x,y) \in D\big\} = \{0\}\\
&\tau^\nu\mbox{-}\nargmin_{y\in Y^\nu}\big\{g^\nu(x,y)~\big|~H^\nu(x,y) \in D\big\} = \emptyset
\end{align*}
regardless of $\tau,\tau^\nu \in [0,\infty)$. Thus, {\rm (NP)}$^\nu$ cannot lead to a lower bound on $\mathfrak{m}$ because its minimum value  is $\infty$ due to the infeasibility of the lower-level problem. In contrast, Theorem \ref{thm:withoutcalmness}(a) applies. 

The example also highlights the critical role of $u$ in {\rm (P)}$^\nu$: If $u$ is removed by fixing it to $0$, then {\rm (P)}$^\nu$ becomes infeasible.   
\end{example}
\noindent {\bf Detail}. Suppose that $f(x,y) = f^\nu(x,y) = x+y$ and $X = [-1,1]$. Then, $\mathfrak{m} = -1$ and $\mathfrak{m}^\nu = -1 - 1/(4\sigma^\nu) + \sigma^\nu/\nu$. We see that even if $\{\sigma^\nu, \nu\in\nats\}$ is bounded, $\nliminf \mathfrak{m}^\nu \leq \mathfrak{m}$. The lower bound improves in the limit to coincide with $\mathfrak{m}$ if $\sigma^\nu\to \infty$ and $\sigma^\nu/\nu\to 0$; cf. Assumption \ref{ass:basic}(f).\eop

Theorem \ref{thm:withoutcalmness}(b) asserts that (P)$^\nu$ may also furnish an upper bound on the minimum value of (P) in the limit. Now,  $\sigma^\nu$ and $\theta^\nu$ tend to infinity, which in view of Assumption \ref{ass:basic}(f) requires coordination with $\delta^\nu$ and $\eta^\nu$; see the discussion after Assumption \ref{ass:basic}. The existence of a compact set $B$ is a nontrivial assumption and practically might only be empirically verified during computations. The main challenge relates to the variable $\lambda$ in (P)$^\nu$, which acts as a penalty parameter in a penalization of $H^\nu(x,z) \in D$; see \eqref{eqn:PnuConstr_mu}. Thus, bounds on $\lambda$ is connected to exact penalization in the lower-level problem. This is trivial if $Y$ is finite as in integer programs with a bounded feasible set. We return to this topic below.

Theorem \ref{thm:withoutcalmness}(c) establishes that a sequence of near-optimal solutions of (P)$^\nu$ can only have cluster points that are feasible for (P), with the tolerance in solving (P)$^\nu$ potentially remaining bounded away from zero.  The existence of a cluster point is not guaranteed by Theorem \ref{thm:withoutcalmness}(c), and it might simply be observed empirically. This is again related to exact penalization in the lower-level problem as we discuss in Subsections \ref{subsec:calmness} and \ref{subsec:localcalmness}. 
In addition to the fact recorded at the end of the theorem for the case $\tau \in (0,\infty)$, one satisfies \eqref{eqn:extraass} regardless of $\tau \in [0,\infty)$ if there is $x\in X$ such that $H(x,z) \in D$ and $f(x,z) < \infty$ for all $z\in Y$ and $\inf_{z\in Y} g(x,z)$ is attained.

\subsection{Bounds and Optimal Solutions under Calmness at a Point}\label{subsec:calmness}

A challenge when relating (P) to (P)$^\nu$ is the fact that the lower-level problem of the former involves the constraint $H(x,z) \in D$, while the latter has a term $\lambda \dist(H^\nu(x,z), D)$. Still, under broad conditions, 
\[
\inf_{z\in Y} \big\{g(x,z) \,\big|\, H(x,z) \in D\big\} = \inf_{z\in Y} g(x,z) + \lambda \dist\big(H(x,z), D\big)~~\mbox{ for sufficiently large } \lambda
\]
according to the principle of exact penalization. Despite the presence of approximations and other changes, one can glean from the inequality constraints \eqref{eqn:PnuConstr_mu} in (P)$^\nu$ a role for exact penalization. This subsection examines the stronger conclusions that become available under a condition for exactness.    

For $x\in \reals^n$ and $u\in \reals^q$, we adopt the notation
\[
V(x,u) = \inf_{y\in Y} \big\{g(x,y)~\big|~H(x,y) + u \in D \big\},
\]
with $V$ representing the {\em value function} for the lower-level problem in (P) under changes to $u$. 

\begin{definition}{\rm (calmness at a point).}\label{def:calmness}
  The value function $V$ is {\em calm at} $x\in \reals^n$ if there is $\lambda \in [0,\infty)$ such that 
  \[
  V(x,u) \geq V(x,0) - \lambda \|u\| ~~~\forall u\in \reals^q,    
  \]
  with $\lambda$ then being referred to as a {\em penalty threshold}.
\end{definition}

We refer to \cite{Burke.91,Burke.91b} for a classical treatment of calmness in constrained optimization, which is closely related to the development here. As compared to \cite{Burke.91,Burke.91b}, we focus on minimum values instead of locally minimum values and we also need to consider the role of $x$. Thus, our definitions and results about exact penalization deviate slightly from those in the literature and proofs are given in Section \ref{sec:proofs}. 

The norm in Definition \ref{def:calmness} is the same as the one underpinning the point-to-set distance. This yields the following characterization of calmness at a point. 

\begin{proposition}{\rm (characterization of calmness at a point).}\label{prop:char_calm}
  For nonempty $D$ and $x\in \reals^n$, the value function $V$ is calm at $x$ with penalty threshold $\lambda$ if and only if 
  \[
  \inf_{y\in Y} g(x,y) + \lambda \dist\big(H(x,y), D\big) = \inf_{y\in Y} \big\{g(x,y) ~\big|~ H(x,y) \in D\big\}. 
  \]
\end{proposition}

In view of the proposition, calmness at $x$ is therefore exactly the property that is needed for exact penalization of the constraint $H(x,z)\in D$ in the lower-level problem of (P) for $x$. We return to a sufficient condition for calmness after the main theorem of this subsection.

\begin{theorem}{\rm (bounds and optimality under calmness).}\label{thm:withcalmness} Suppose that $\theta^\nu|\tau^\nu - \tau| \to 0$ and there exists an optimal solution $(x^\star,y^\star)$ of {\rm (P)} such that the value function $V$ is calm at $x^\star$. Under Assumption \ref{ass:basic}, the following hold: 
\begin{enumerate}[(a)]

\item If $\{\theta^\nu, \nu\in\nats\}$ is bounded away from zero, then $\nlimsup \mathfrak{m}^\nu \leq \mathfrak{m} \in \reals$. 

\item Suppose that $\sigma^\nu, \theta^\nu\to \infty$ and there is a compact set $B\subset \reals^{n+m+q+2}$ such that for each $\nu\in\nats$ there exists an optimal solution of {\rm (P)}$^\nu$ in $B$. Let $\lambda^\star \in [0,\infty)$ be a penalty threshold for the value function $V$ at $x^\star$. If $\lambda\in [\lambda^\star,\infty)$, then there exist a vanishing sequence $\{\epsilon^\nu\in [0,\infty), \nu\in\nats\}$ and a sequence $\{w^\nu\in\reals^{n+m+q+2}, \nu\in\nats\}$ such that $w^\nu$ is an $\epsilon^\nu$-optimal solution of {\rm (P)}$^\nu$ and 
    \[
    w^\nu \to (x^\star, y^\star, 0, 0, \lambda).
    \]
    Moreover,  $\nlim \mathfrak{m}^\nu = \mathfrak{m} \in \reals$. 

\item Suppose that $\sigma^\nu\to\infty$, $\theta^\nu\to \infty$, $\epsilon^\nu \to 0$, and, for each $\nu\in\nats$, $(x^\nu, y^\nu, u^\nu, \alpha^\nu, \lambda^\nu)$ is an $\epsilon^\nu$-optimal solution of {\rm (P)}$^\nu$. If there is a subsequence $N\subset\nats$ and $(\hat x, \hat y, \hat u, \hat \alpha,$ $\hat \lambda)$ such that 
    \[
    (x^\nu, y^\nu, u^\nu, \alpha^\nu, \lambda^\nu)\Nto (\hat x, \hat y, \hat u, \hat \alpha, \hat \lambda),
     \]
     then $(\hat x, \hat y)$ is an optimal solution of {\rm (P)} and $\mathfrak{m}^\nu \Nto \mathfrak{m} \in \reals$. 

\end{enumerate}

\end{theorem}

Theorem \ref{thm:withcalmness}(a) supplements Theorem \ref{thm:withoutcalmness}(a) by not requiring a compact $Y$ and allowing $\tau^\nu = \tau$ for all $\nu$. The requirement about calmness is mild as it only needs to hold at {\em some} $x^\star$, and this point does not need to be known when solving (P)$^\nu$. 
Again, an examination of the proof of Theorem \ref{thm:withcalmness}(a) shows that the liminf-condition in Assumption \ref{ass:basic}(e) is superfluous; cf. the discussion after Theorem \ref{thm:withoutcalmness}.

\begin{example}{\rm (lack of calmness).}\label{ex:calmness}
For the data in Example \ref{ex:infeaslower}, $V(x,u) = 0$ if $u\leq 0$ and $V(x,u) = \infty$ otherwise. Since the value function $V$ in this case is independent of $x$, it is calm at every $x$. The situation changes if $g(x,y) = y$. Then, 
\[
V(x,u) = \begin{cases}
  \infty & \mbox{ if } u >0\\
  -\sqrt{-u} & \mbox{ if } u \in [-1,0]\\
  -1 & \mbox{ otherwise,}
\end{cases}
\]
which is not calm at any $x$. This might invalidate the conclusion $\nlimsup \mathfrak{m}^\nu \leq \mathfrak{m}$ in Theorem \ref{thm:withcalmness}(a). 
\end{example}
\noindent {\bf Detail}. Suppose that $g^\nu(x,y) = g(x,y) = y$ so that $\delta^\nu =0$ in Assumption \ref{ass:basic} and $\theta^\nu = \nu^{1/3}$, $\sigma^\nu = \nu^{2/3}$, $\bar \lambda^\nu = \nu^{1/3}$, which implies that Assumption \ref{ass:basic}(f) holds because $\eta^\nu = 1/\nu$. In this case with $\tau^\nu = 0$, (P)$^\nu$ takes the form 
\begin{align*}
  &\nnmin_{x,y,u,\alpha,\lambda} \, f^\nu(x,y) + \nu^{2/3} |u| - \nu^{1/3}\alpha\\
  &\mbox{ subject to }   y^2 + 1/\nu + u \leq 0, ~x\in X, ~y\in [-1,1], ~\lambda \in [0,\nu^{1/3}], ~\alpha \leq 0\\
  &~~~~~~~~~~~~~~~ y + \alpha \leq z + \lambda \max\{0, z^2 + 1/\nu\} ~~~\forall z\in [-1,1].
\end{align*}  
For $\lambda \geq 1/2$, the last constraints reduce to the single constraint $y + \alpha \leq -1/(4\lambda) + \lambda/\nu$. Since $\lambda = \bar \lambda^\nu$ is always optimal, the constraint simplifies further to 
\[
y + \alpha \leq -\frac{1}{4\nu^{1/3}} + \frac{1}{\nu^{2/3}}. 
\]
While the right-hand vanishes, it has a negative value of order $\nu^{-1/3}$ for large $\nu$. The constraint can be satisfied by setting $y=0$, but then $\alpha$ is of order $\nu^{-1/3}$, which incurs a nonvanishing penalty in the objective function because  $\theta^\nu = \nu^{1/3}$. Alternatively, $\alpha =0$ and $y$ is of order $\nu^{-1/3}$ but then $u$ must be of order $\nu^{-2/3}$. Since $\sigma^\nu = \nu^{2/3}$, we again face a nonvanishing penalty term in the objective function. Consequently, even if $f^\nu$ is identical to $f$, $\mathfrak{m}$ remains strictly below $\mathfrak{m}^\nu$ and the conclusion of Theorem \ref{thm:withcalmness}(a) fails.\eop

Theorem \ref{thm:withcalmness}(b) has no counterpart in Theorem \ref{thm:withoutcalmness}. It establishes that {\em any} optimal solution of (P) satisfying the calmness assumption of the theorem can be recovered by solving (P)$^\nu$. Thus, there is no optimal solution of (P) of this kind that is ``isolated'' in the sense that it is not approachable by near-optimal solutions of (P)$^\nu$. The $\lambda$-components of the stipulated sequence $\{w^\nu, \nu\in\nats\}$ are bounded. This partially addresses the concern from the discussion after  Theorem \ref{thm:withoutcalmness} about unbounded $\lambda$-values. 

The conclusion in Theorem \ref{thm:withcalmness}(c) strengthens that of Theorem \ref{thm:withoutcalmness}(c) by guaranteeing optimality in (P), and not merely feasibility. In addition, the minimum values $\mathfrak{m}^\nu$ of (P)$^\nu$ now converge, at least along the subsequence indexed by $N$, to the minimum value of (P).  

As an examples of a sufficient condition for calmness at a point, we furnish the following proposition. We denote by  
$\dist_2(y,C) = \inf_{y'\in C} \|y - y'\|_2$ the {\em point-to-set distance in the Euclidean norm}. 

\begin{proposition}{\rm (sufficient condition for calmness at a point).}\label{prop:suffCalm} For the value function $V$ to be calm at $x\in \reals^n$, the following suffices: There exist $\kappa_1,\kappa_2 \in [0,\infty)$ such that 
\begin{align}
&\big|g(x,y) - g(x,y')\big| \leq \kappa_1 \|y - y'\|_2 ~~\forall y,y'\in Y\label{eqn:LipAssump}\\
& \dist_2\big(y, S(x)\big) \leq \kappa_2 \dist\big(H(x,y), D\big) ~~\forall y\in Y,\label{eqn:regularity}\\
&~~~~~~~~~~~ \mbox{ where } S(x) = \big\{z\in Y~\big|~H(x,z) \in D\big\}.\nonumber
\end{align}
\end{proposition}

The regularity condition \eqref{eqn:regularity} specifies a certain growth rate of the constraint mapping $H(x, \cdot)$ as $y$ moves away from the lower-level feasible set $S(x)$. This holds trivially if $Y$ is finite as in some integer programs. It also holds when $Y$ is a polyhedron, $H(x,\cdot)$ is affine, and $D$ is polyhedral. Then, Hoffman's lemma \cite{Hoffman.52} ensures the existence of $\kappa_2$ such that \eqref{eqn:regularity} holds; see, for example, \cite{PenaVeraZuluaga.21} for recent advances in quantifying $\kappa_2$, which is called the Hoffman's constant in this case. More generally, \eqref{eqn:regularity} relates to constraint qualifications as discussed in the next subsection. 

Growth of the values of a constraint mapping at points moving away from the feasible set $S(x)$ is distinct from growth of lower-level objective function values at points moving away from lower-level optimal solutions. The latter is discussed in \cite{YeZhu.95} under the name uniformly weak sharp minima.

\subsection{Optimal Solutions under Local Calmness}\label{subsec:localcalmness}

We obtain an even stronger conclusion if a {\em local} calmness property holds at a cluster point of solutions generated from (P)$^\nu$. This largely eliminates concerns about unbounded $\lambda$-values.

\begin{definition}{\rm (local calmness).}\label{def:localcalm}
  The value function $V$ is {\em locally calm at} $\bar x\in \reals^n$ if there is $\lambda \in [0,\infty)$ and $\rho\in (0,\infty)$ such that 
  \[
  V(x,u) \geq V(x,0) - \lambda \|u\| ~~~\forall u\in \reals^q, ~x\in \ball(\bar x, \rho),   
  \]
  with $\lambda$ being referred to as a {\em penalty threshold} and $\rho$ as a {\em radius of calmness}.
\end{definition}

Local calmness at $\bar x$ would hold, for example, if the constants $\kappa_1,\kappa_2$ in Proposition \ref{prop:suffCalm} apply not only at a point, but also locally near $\bar x$. Again, in the polyhedral case, Hoffman's lemma \cite{Hoffman.52} enters with the uniformity locally in $x$ corresponding to having an upper bound on the Hoffman's constants as the polyhedral system is perturbed by $x$. This is a relatively mild assumption because, for example, the Hoffman's constant of a linear system of inequalities is independent of right-hand side perturbations. 

In the case of $D = (-\infty,0]^q$, $Y = \{y\in \reals^m~|~g_i(y) \leq 0, i=1, \dots, r\}$ being compact for smooth $g_i:\reals^m\to \reals$, $\{y\in Y~|~H(x,y) \in D\} \neq \emptyset$ for all $x\in \reals^n$, and $H$ is smooth, Theorem 2.10 in \cite{RoysetPolakDerkiureghian.03} establishes that the value function $V$ is locally calm at every $\bar x\in \reals^n$ under the constraint qualification: for all $x\in \reals^n$ and $y\in \nargmin_{z\in Y} \{g(x,z)~|~H(x,z) \in D\}$, one has 
\[
  \big\{\nabla g_i(y), i \in I(y); ~~\nabla_y h_j(x,y), j\in J(x,y)\big\} ~~\mbox{ are linearly independent vectors,}
\]
where $I(y) = \{i~|~ g_i(y) = 0\}$, $J(x,y) = \{j~|~h_j(x,y) = 0\}$, and $h_j$ is the $j$th component of $H$.  

We refer to \cite{Burke.91} and references therein for other constraint qualifications related to local calmness.

\begin{theorem}{\rm (optimality under local calmness).}\label{thm:withLocalCalmness} Suppose that $\sigma^\nu,\theta^\nu\to \infty$ while $\epsilon^\nu$ and $\theta^\nu|\tau^\nu - \tau|$ vanish, Assumption \ref{ass:basic} holds, and there exists an optimal solution $(x^\star,y^\star)$ of {\rm (P)} such that the value function $V$ is calm at $x^\star$. Moreover, for each $\nu\in\nats$, let $(x^\nu, y^\nu, u^\nu, \alpha^\nu, \lambda^\nu)$ be an $\epsilon^\nu$-optimal solution of {\rm (P)}$^\nu$. Suppose there is a subsequence $N\subset\nats$ and $(\hat x, \hat y, \hat u, \hat \alpha)$ such that 
\[
(x^\nu, y^\nu, u^\nu, \alpha^\nu)\Nto (\hat x, \hat y, \hat u, \hat \alpha)
\]
and the value function $V$ is locally calm at $\hat x$. Then, the following hold: 
\begin{enumerate}[(a)]

\item If $Y^\nu = Y$ for all $\nu$ sufficiently large, then $(\hat x, \hat y)$ is an optimal solution of {\rm (P)} and $\mathfrak{m}^\nu \Nto \mathfrak{m}\in\reals$. 

\item If $\theta^\nu\bar \lambda^\nu \sup_{y\in Y} \dist_2(y,Y^\nu) \to 0$ and there are Lipschitz constants $\kappa_g, \kappa_H\in [0,\infty)$ such that 
\begin{align*}
  \big| g(x,y) - g(x,y') \big| & \leq \kappa_g \|y - y'\|_2  ~~~\forall x\in X, ~y,y'\in Y\\ 
  \big\| H(x,y) - H(x,y') \big\| & \leq \kappa_H \|y - y'\|_2  ~~~\forall x\in X, ~y,y'\in Y, 
\end{align*}

then $(\hat x, \hat y)$ is an optimal solution of {\rm (P)} and $\mathfrak{m}^\nu \Nto \mathfrak{m}\in \reals$.

\end{enumerate}

\end{theorem}

Theorem \ref{thm:withLocalCalmness} is distinct from Theorems \ref{thm:withoutcalmness}(c) and \ref{thm:withcalmness}(c) by {\em not} requiring that the $\lambda$-components of the near-optimal solutions of (P)$^\nu$ are bounded, at least along a subsequence. It is only the $x$-, $y$-, $u$-, and $\alpha$-components that need to have a cluster point in Theorem \ref{thm:withLocalCalmness}. This is nearly automatic, for example, when $X$ and $Y$ are compact. Local calmness is only required to hold at the $x$-component of such a cluster point, while pointwise calmness is assumed to hold at {\em some} $x^\star$.   

When $Y^\nu = Y$, no other assumptions are needed to ensure that the $x$- and $y$-components of a cluster point form an optimal solution of (P); see part (a). When $Y^\nu$, $\nu\in\nats$, are strict subsets of $Y$, additional Lipschitz properties enter in (b). While these might be perceived as mild, we also rely on a sufficiently slow growth of $\theta^\nu$ and $\bar \lambda^\nu$ relative to the convergence of $Y^\nu$ to $Y$, which could be more restrictive.

\section{Illustrative Examples}\label{sec:example}

Example \ref{ex:infeaslower} illustrates the importance of introducing the additional vector $u$ in (P)$^\nu$. The use of $\lambda$ for penalization of lower-level constraints in (P)$^\nu$ is computationally motivated: Without penalization the index set of the constraints \eqref{eqn:PnuConstr_mu} would become $x$-dependent resulting in a challenging {\em generalized} semi-infinite problem. The importance of the final addition in (P)$^\nu$ relative to (P), the variable $\alpha$, emerges from the next example. 

\begin{example}{\rm (role of $\alpha$).}\label{ex:alphaimp}
For $g^\nu(x,y) = e^y/\nu$, $Y=Y^\nu = (-\infty,0]$, $H^\nu(x,y) = 0$, and $D = \{0\}$, we obtain 
\[
\inf_{y\in Y^\nu} g^\nu(x,y) + \lambda \dist\big(H^\nu(x,y), D\big) = 0
\]
for any $\lambda \in [0,\infty)$. Thus, there is no $y\in Y$ satisfying the constraints \eqref{eqn:PnuConstr_mu} when $\alpha = 0$. The introduction of this variable therefore allows for satisfaction of \eqref{eqn:PnuConstr_mu} without imposing additional assumptions.
\end{example}

We next consider a more fully developed example where arbitrarily small changes to (P) cause large changes in its solution while (P)$^\nu$ is better behaved.  
\begin{example}{\rm (instability in (P)).}\label{ex:instab}
Let $X = [1,2]$, $Y = \{0,1\}$, $D = (-\infty,0]$, $f(x,y) = (y-1/2)x$, $g(x,y) = -xy$, and $H(x,y) = y - 1$. Since $\nargmin_{y\in Y}\{g(x,y)~|~$ $H(x,y) \in D\} = \{1\}$ for any $x\in X$, (P) under tolerance $\tau = 0$ reduces to 
minimizing $(y-1/2)x$ over $x\in [1,2]$ and $y\in \{0,1\}$ while satisfying $y = 1$. Thus, the unique optimal solution of {\rm (P)} is $(x^\star, y^\star) = (1,1)$ with minimum value $\mathfrak{m} = 1/2$. 

For $\nu\in\nats$, suppose that $H(x,y)$ is approximated by $H^\nu(x,y)  = y - 1 + 1/\nu$, while the other components remain unchanged. Thus, parts (a,c,d,e) of Assumption \ref{ass:basic} hold with $\delta^\nu = 0$ and $\eta^\nu = 1/\nu$. The naive approach of simply substituting in the approximating quantities produces the following instance of {\rm (NP)}$^\nu$:   
\begin{align*}
\nnmin_{x\in X, y\in \reals}  ~ f^\nu(x,y)=(y-\tfrac{1}{2})x ~\mbox{ subject to } ~ &y \in \nargmin_{z\in Y} \big\{ g^\nu(x,z) ~\big|~ H^\nu(x,z) \in D\big\}\\
&=\nargmin_{z\in \{0,1\}} \{ -xz ~|~ z - 1 + 1/\nu \leq 0\}, 
\end{align*}
with unique optimal solution $(\tilde x^\nu, \tilde y^\nu) = (2,0)$ and corresponding minimum value $\tilde{\mathfrak{m}}^\nu = -1$. As $\nu\to \infty$, neither the optimal solutions nor the minimum values of {\rm (NP)}$^\nu$ converge to the corresponding quantities of {\rm (P)}. This occurs despite a seemingly minor discrepancy in $H^\nu$ relative to $H$. It turns out that {\rm (P)}$^\nu$ is much better behaved. 
\end{example}
\noindent {\bf Detail}. For $Y^\nu = Y$ and $\tau^\nu = 0$, it takes the form 
\begin{align*}
\nnmin_{x,y,u,\alpha,\lambda}  ~ (y-\tfrac{1}{2})x + \sigma^\nu |u| & - \theta^\nu \alpha\\
 \mbox{ subject to } ~ y - 1 + 1/\nu + u & \leq 0\\
 - xy + \alpha & \leq - xz + \lambda \max\{0, z - 1 + 1/\nu\} ~~~\forall z\in \{0,1\}\\
 x\in [1,2], ~y&\in \{0,1\}, ~u\in \reals, ~\alpha  \leq 0, ~\lambda \in [0, \bar\lambda^\nu].
\end{align*}
In light of Assumption \ref{ass:basic} and the theorems above, the parameters $\sigma^\nu$, $\theta^\nu$, and $\bar\lambda^\nu$ can be set to 
$\sigma^\nu = \nu^{1/2}$, $\theta^\nu = \nu^{1/3}$, and $\bar\lambda^\nu = \nu^{1/3}$. Since $\delta^\nu = 0$ and $\eta^\nu = 1/\nu$, Assumption \ref{ass:basic}(f) holds while both $\sigma^\nu$ and $\theta^\nu$ tend to infinity as $\nu\to \infty$. An optimal solution of this instance of (P)$^\nu$ is determined, for example, by recognizing that $\lambda$ can be set to $\bar\lambda^\nu= \nu^{1/3}$ and the two constraints indexed by $z$ reduces to the single constraint $-xy + \alpha \leq -x + \nu^{-2/3}$. This leads to the optimal solutions and minimum values: 
\begin{align*}
&\nu = 1, \dots, 6\hspace{-0.1cm}: (x^\nu,y^\nu,u^\nu,\alpha^\nu,\lambda^\nu) = (1, 0, 0, -1 + \nu^{-2/3}, \nu^{1/3}), \mathfrak{m}^\nu = -\tfrac{1}{2} + \nu^{1/3} - \nu^{-1/3}\\
&\nu = 7, 8, \dots: (x^\nu,y^\nu,u^\nu,\alpha^\nu,\lambda^\nu) = (1, 1, -1/\nu, 0, \nu^{1/3}), ~~~~\,~~~\mathfrak{m}^\nu = \tfrac{1}{2} + \nu^{-1/2}.
\end{align*} 
Clearly, $(x^\nu,y^\nu)\to (x^\star, y^\star)$ and $\mathfrak{m}^\nu\to \mathfrak{m}$. 

In this instance of (P), the value function $V$ is actually locally calm at any $x\in \reals$. Thus, there is no need for the $\lambda$-component of a solution of (P)$^\nu$ to grow beyond bounds and, in fact, it can remain at zero. This insight produces another set 
of optimal solutions and minimum values of (P)$^\nu$:
\begin{align*}
&\mbox{for } \nu = 1, 2, 3 ~~~~~~\,(x^\nu,y^\nu,u^\nu,\alpha^\nu,\lambda^\nu) = (1, 0, 0, -1, 0), ~~~~~~~~~\mathfrak{m}^\nu = -\tfrac{1}{2} + \nu^{1/3}\\
&\mbox{for } \nu = 4, 5, \dots  ~~~~(x^\nu,y^\nu,u^\nu,\alpha^\nu,\lambda^\nu) = (1, 1, -1/\nu, 0, 0), ~~~~~~\mathfrak{m}^\nu = \tfrac{1}{2} + \nu^{-1/2},
\end{align*} 
with slightly faster recovery of $(x^\star, y^\star)$ as $\nu\to \infty$.\eop 

The previous example involves a nonconvex set $Y = \{0,1\}$. However, difficulties are not limited to nonconvex lower-level problems as the next example shows.

\begin{example}{\rm (instability of (P) under convex $Y$).}\label{ex:instab2}
For $g(x,y) = 0$, $g^\nu(x,y) = y^2/\nu$, $Y=Y^\nu = [-2,1]$, $H(x,y) = H^\nu(x,y) = -y-1$, $D = (-\infty, 0]$, $X = [0, 1]$, and $f(x,y) = f^\nu(x,y) = (y+1/2)x$, we obtain 
\begin{align*}
&\nargmin_{y\in Y} \big\{ g(x,y) ~\big|~ H(x,y) \in D \big\} = [-1,1]\\
&\nargmin_{y\in Y} \big\{ g^\nu(x,y) ~\big|~ H^\nu(x,y) \in D \big\} = \{0\},
\end{align*}
with these lower-level problems being convex. Suppose that $\tau = \tau^\nu = 0$. Thus, we find that {\rm (P)} has minimum value $\mathfrak{m} = -1/2$, with optimal solution $x^\star = 1$, $y^\star = -1$. The corresponding numbers for the naive problem {\rm (NP)}$^\nu$ are 0, 0, and 0, respectively, regardless of $\nu$. 

Again, {\rm (P)}$^\nu$ is better behaved. If $\theta^\nu = \sqrt{\nu}$, then Theorem \ref{thm:withLocalCalmness} applies with calmness and local calmness holding trivially.  
\end{example}
\noindent {\bf Detail}. For $\sigma^\nu > 1$, every optimal solution of (P)$^\nu$ has $u = 0$. In \eqref{eqn:PnuConstr_mu}, $z = 0$ is an option which makes the right-hand side of those constraint equal to zero. Thus, (P)$^\nu$ reduces to 
\[
\nnmin_{x,y,\alpha} ~(y+1/2)x - \theta^\nu \alpha ~\mbox{ subject to } ~x \in [0, 1], ~y\in [-1, 1], ~y^2/\nu + \alpha \leq 0, ~\alpha \leq 0. 
\]  
By inspection, $(0,0,0)$ is an optimal solution for $\nu = 1, 2, 3$. For $\nu\geq 4$, this changes to $x^\nu =1$, $y^\nu = -1$, $\alpha^\nu = -1/\nu$, with minimum value $\mathfrak{m}^k = -1/2 + 1/\sqrt{\nu}$.\eop

We refer to \cite[Section 4.4.2]{CarusoLignolaMorgan.20} for other examples.

\section{Computational Considerations}\label{sec:computational}

Next, we fix $\nu\in\nats$ and examine the challenges associated with solving (P)$^\nu$. This provides the basis for an algorithm for (P) that iteratively and approximately solves (P)$^\nu$ for increasingly large $\nu$. We start by establishing sufficient conditions for the existence of solutions before turning to an outer approximation algorithm. 

\begin{proposition}{\rm (existence of near-optimal solutions for (P)$^\nu$).} 
For $f^\nu:\reals^n\times\reals^m\to (-\infty,\infty]$, $g^\nu: \reals^n\times\reals^m\to \reals$, and $H^\nu:\reals^n\times\reals^m\to \reals^q$, suppose that $D$ and $X \times Y  \cap \dom f^\nu$ are nonempty and, for all $x\in X$, 
\begin{equation}\label{eqn:existence}
\inf_{z\in Y^\nu} g^\nu(x,z) + \bar \lambda^\nu \dist\big(H^\nu(x,z), D\big) > - \infty. 
\end{equation}
Then, the following hold: 
\begin{enumerate}[(a)]

\item If $\epsilon^\nu \in (0,\infty)$ and $\inf_{x\in X, y\in Y} f^\nu(x,y) > -\infty$, then there exists an $\epsilon^\nu$-optimal solution of {\rm (P)}$^\nu$. 

\item If $\sigma^\nu>0$, $\theta^\nu>0$, the functions $f^\nu, g^\nu$ are lsc, $H^\nu$ is continuous relative to $X\times Y$, the function $g^\nu(\cdot,z)$ is continuous relative to $X$ for each $z\in Y^\nu$, the sets $X, Y, D$ are closed, and the set
\begin{equation}\label{eqn:levelsetass}
\big\{(x,y) \in X \times Y~\big|~f^\nu(x,y) \leq \gamma\big\}
\end{equation}
is bounded for each $\gamma\in \reals$, then there exists an optimal solution of {\rm (P)}$^\nu$. 

\end{enumerate}
\end{proposition}
\noindent {\bf Proof}. For (a), it suffices to confirm that $\mathfrak{m}^\nu$ is finite. By assumption, there exist $\bar x\in X$, $\bar y \in Y$, $\bar u\in \reals^q$, and $\bar \alpha\in (-\infty,0]$ such that $f^\nu(\bar x, \bar y) < \infty$, $H^\nu(\bar x, \bar y) + \bar u \in D$, and  
\[
g^\nu(\bar x, \bar y) + \bar\alpha \leq \inf_{z\in Y^\nu} g^\nu(\bar x,z) + \bar \lambda^\nu \dist\big(H^\nu(\bar x,z), D\big) + \tau^\nu.
\]
Thus, $(\bar x, \bar y, \bar u, \bar \alpha, \bar\lambda^\nu)$ is feasible in (P)$^\nu$ and $\mathfrak{m}^\nu < \infty$. Since $\inf_{x\in X, y\in Y} f^\nu(x,y) > -\infty$ and the last two terms in the objective function of (P)$^\nu$ are nonnegative at every feasible point, $\mathfrak{m}^\nu >- \infty$. 

For (b), we have $\mathfrak{m}^\nu < \infty$ by following the arguments in (a). The objective function in (P)$^\nu$ is lsc and level-bounded on the feasible set because $\sigma^\nu>0$, $\theta^\nu>0$, and \eqref{eqn:levelsetass} is bounded. By \cite[Theorem 4.9]{primer} it only remains to show that the feasible set is closed. Since an intersection of closed sets is closed, this also holds straightforwardly; we recall that $u\mapsto \dist(u,D)$ is continuous.\eop

The assumptions in the proposition are mild. For example, if $g^\nu(x,\cdot)$ is lsc and $Y^\nu$ is compact, then \eqref{eqn:existence} holds. A lower bound on $f^\nu$ on $X\times Y$ might also be readily available. Thus, a near-optimal solution of (P)$^\nu$ can be expected via (a) and this is all what is needed in the earlier theorems. For general existence results in bilevel optimization, we refer to \cite{CarusoLignolaMorgan.20}.

While the existence of near-optimal solutions of (P)$^\nu$ is reassuring, we are still faced with the nontrivial task of computing them. Naturally, much depends on any additional properties the components of the problem might have and a comprehensive treatment extends beyond the present paper. We provide a generic algorithm based on outer approximations that, in part, addresses the potentially large, even infinite, number of constraints in \eqref{eqn:PnuConstr_mu}. We adopt the following notation to describe (P)$^\nu$: 
\begin{align*}
&C = \big\{ (x,y,u,\alpha,\lambda) \in \reals^{n+m+q+2}~\big|~x\in X, y\in Y, \alpha \leq 0, \lambda \in [0,\bar \lambda^\nu],\\
&~~~~~~~~~~~~~~~~~~~~~~~~~~~~~~~~~~~~~~~~~~~~~ H^\nu(x,y) + u \in D, (x,y) \in \dom f^\nu\big\}\\
&\psi_0(x,y,u,\alpha,\lambda) = f^\nu(x,y) + \sigma^\nu \|u\| - \theta^\nu \alpha\\
&\psi_z(x,y,u,\alpha,\lambda) = g^\nu(x,y) + \alpha - g^\nu(x,z) - \lambda \dist\big(H^\nu(x,z), D\big) - \tau^\nu, ~~~~z\in Y^\nu.
\end{align*}
Thus, (P)$^\nu$ amounts to minimizing $\psi_0(w)$ subject $w\in C$ and $\psi_z(w) \leq 0$ for $z\in Y^\nu$. Since $\nu$ is fixed in this section, we suppress the dependence on $\nu$ in $\psi_0, \psi_z$, and $C$.

\medskip

\noindent  Outer Approximation Algorithm for (P)$^\nu$.

\begin{description}

  \item[Data.] ~\,~$(x^0, y^0, u^0, \alpha^0, \lambda^0) \in C$, $Y^\nu_0 = \emptyset$, $\{\epsilon^k, \delta^k \geq 0, k\in \nats\}$.

  \item[Step 0.]  Set $k = 1$.

  \item[Step 1.]  Set $Y^\nu_k = Y^\nu_{k-1} \cup \{z^k\}$, where 
  \[
  z^k \in \delta^k\mbox{-}\nargmin_{z\in Y^\nu} g^\nu(x^{k-1},z) + \lambda^{k-1} \dist\big(H^\nu(x^{k-1},z), D\big).
  \] 

  \item[Step 2.]  Compute $(x^k, y^k, u^k, \alpha^k, \lambda^k) \in \epsilon^k\mbox{-}\nargmin_{w \in C} \big\{\psi_0(w)~\big|~\psi_z(w) \leq 0, z\in Y_k^\nu\big\}$.

\item[Step 3.] Replace $k$ by $k +1$ and go to Step 1.
\end{description}

The initialization of the algorithm is usually straightforward as any $(x,y)\in (X\times Y) \cap \dom f^\nu$ can be paired with some $u$, which is a vector of slack-type variables, and an arbitrary $\alpha\leq 0$ and $\lambda \in [0,\bar\lambda^\nu]$ to produce a point in $C$. 
Step 1 identifies a constraint from \eqref{eqn:PnuConstr_mu} to be added to the collection of constraints indexed by $Y^\nu_{k-1}$. If $D = (-\infty, 0]^q$ and $\|\cdot\| = \|\cdot\|_1$, then $\dist(H^\nu(x^{k-1},z), D)$ equals 
\begin{equation}\label{eqn:onenormexam}
\sum_{j = 1}^q \max\big\{0, h^\nu_j(x^{k-1},z)\big\},
\end{equation}
where $h_j^\nu$ is the $j$th component of $H^\nu$. Moreover, if $g^\nu(x^{k-1},\cdot)$ and $h_j^\nu(x^{k-1}, \cdot)$, $j=1, \dots, q$, are affine functions and $Y^\nu$ is polyhedral, possibly restricted further to integer variables, then the resulting subproblem can be solved using (mixed-integer) linear programming. 

Step 2 amounts to solving a master problem similar to (P)$^\nu$, but with only $k$ constraints of the form \eqref{eqn:PnuConstr_mu}. Again, we permit an optimality tolerance. A potential challenge is the variable $\lambda$, which is multiplied by $\dist(H^\nu(x,z),D)$ and thus produces a nonlinearity. However, we know that $\psi_z(x,y,u,\alpha,\lambda) \leq \psi_z(x,y,u,\alpha,\lambda')$ if $\lambda \geq \lambda'\geq 0$, which implies that one can fix $\lambda = \bar\lambda^\nu$ in Step 2. (A reason for {\em not} fixing $\lambda$ at its upper limit would be numerical instability issues, which might trigger substeps  fixing $\lambda$ at a lower level initially and then increase it cautiously.) 

A more serious challenge is the term $-\dist(H^\nu(x,z),D)$ on the left-hand side in the constraint $\psi_z(w) \leq 0$. In the example leading to \eqref{eqn:onenormexam}, this produces a max-expression on the ``wrong'' side of the inequality. However, such terms can be linearized. {\blue For additional binary variables $w_1, \dots, w_q$ and additional nonpositive continuous variables $v_1, \dots, v_q$, we find that 
\begin{align}
\psi_z(x,y,u,\alpha,\lambda) \leq 0 ~~\Longleftrightarrow~~& g^\nu(x,y) + \alpha - g^\nu(x,z) + \lambda \sum_{j=1}^q v_j  \leq \tau^\nu, ~-U_j w_j \leq v_j,\nonumber\\
& L_j (1-w_j) - h_j^\nu(x,z) \leq v_j \leq -h_j^\nu(x,z)~\forall j\label{eqn:reformconstrs}
\end{align}
provided that $L_j \leq h_j^\nu(x,z) \leq U_j$ for $j=1, \dots, q$.} While there are nontrivial computational challenges associated with the resulting big-M formulation and the growing number of variables and constraints, linearization of this kind provides a concrete path forward to solving the master problem in Step 2 using mixed-integer programming. 

We observe that $\psi_0(x^k, y^k, u^k, \alpha^k, \lambda^k) - \epsilon^k \leq \mathfrak{m}^\nu$ for any $k$, which can be used to bound the suboptimality of any feasible point for (P)$^\nu$. Feasible points are readily available because a slight adjustment of $\alpha^k$ usually leads to a point satisfying \eqref{eqn:PnuConstr_mu}. Regardless, as we see next, the algorithm is guaranteed to achieve an optimal solution of (P)$^\nu$ in the limit under mild assumptions.

\begin{proposition}{\rm (outer approximation algorithm).} Suppose that $X,Y,D$ are nonempty closed sets, $Y^\nu$ is nonempty and compact, 
$g^\nu:\reals^n\times\reals^m\to \reals$ and $H^\nu:\reals^n\times\reals^m\to \reals^q$ are continuous, $f^\nu:\reals^n\times\reals^m\to (-\infty, \infty]$ is continuous relative to $\dom f^\nu$, a closed set, and for each $\bar x\in X$ and $\epsilon \in (0,\infty)$, there is $\rho\in (0,\infty)$ such that for all $z\in Y^\nu$ and $x\in X \cap \ball(\bar x,\rho)$,
\begin{align}
\big|g^\nu(x,z) - g^\nu(\bar x, z)\big| &\leq \epsilon\label{eqn:outerapprox1}\\
\big\|H^\nu(x,z) - H^\nu(\bar x, z)\big\| &\leq \epsilon.\label{eqn:outerapprox2}
\end{align} 
If the outer approximation algorithm generates $(x^k, y^k, u^k, \alpha^k, \lambda^k)$ using vanishing $\epsilon^k$ and $\delta^k$, and 
\[
(x^k, y^k, u^k, \alpha^k, \lambda^k) \Nto (\hat x, \hat y, \hat u, \hat \alpha, \hat \lambda)
\]
for some subsequence $N\subset\nats$, then $(\hat x, \hat y, \hat u, \hat \alpha, \hat \lambda)$ is an optimal solution for {\rm (P)}$^\nu$ and 
\[
\psi_0(x^k, y^k, u^k, \alpha^k, \lambda^k)\Nto \mathfrak{m}^\nu.
\] 
\end{proposition} 
\noindent {\bf Proof}.  The conclusion follows from convergence statement 6.15 in \cite{primer} after verifying the conditions: 
\begin{enumerate}[(i)]
\item $C$ is closed.

\item $\psi_0$ is continuous relative to $C$.

\item $\{\psi_z, z\in Y^\nu\}$ is locally bounded, i.e., for any $\bar w\in \reals^{n+m+q+2}$, there are $\epsilon,\rho\in (0,\infty)$ such that $|\psi_z(w)|\leq \rho$ for all $z\in Y^\nu$ and $w\in \ball(\bar w, \epsilon)$. 

\item $\{\psi_z, z\in Y^\nu\}$ is equicontinuous relative to $C$, i.e., for any $\bar w\in C$ and $\epsilon > 0$ there exists $\rho>0$ such that $|\psi_z(w) - \psi_z(\bar w)|\leq \epsilon$ for all $z\in Y^\nu$ and $w\in C\cap\ball(\bar w, \rho)$. 

\item For any $w\in C$, there is a maximizer of $\psi_z(w)$ over $z\in Y^\nu$. 

\end{enumerate}

Conditions (i) and (ii) hold trivially under the assumed conditions. (We note that \cite[statement 6.15]{primer} requires $\psi_0$ to be real-valued and continuous, but this is easily overcome by redefining $\psi_0$ outside $\dom f^\nu$ using a continuous extension because $\dom f^\nu$ is closed.) 

Condition (iii) holds because $(z,w)\mapsto \psi_z(w)$ is continuous and $Y^\nu$ is compact. 

For condition (iv), let $\bar w = (\bar x, \bar y, \bar u, \bar \alpha, \bar \lambda) \in C$ and $\epsilon > 0$. For any $w=(x,y,u,\alpha,\lambda)$ $\in$ $C$, one has 
\begin{align}\label{eqn:psizexpan}
  &\big|\psi_z(w) - \psi_z(\bar w)\big|  \leq \big|g^\nu(x,y) - g^\nu(\bar x, \bar y)\big| + |\alpha - \bar \alpha| +  \big|g^\nu(x,z) - g^\nu(\bar x, z)\big|\\ 
  & ~~+ \bar\lambda^\nu\Big|\dist\big(H^\nu(x,z),D\big) - \dist\big(H^\nu(\bar x, z),D\big)  \Big| + |\lambda - \bar\lambda | \dist\big(H^\nu(\bar x, z),D\big).\nonumber
\end{align}
Since $H^\nu$ is continuous, $D$ is nonempty, and $Y^\nu$ is compact, one has that $\kappa$ $=$ $\sup_{z'\in Y^\nu} \dist(H^\nu(\bar x, z'),D)$ is finite. Moreover, by continuity of $g^\nu$ and $H^\nu$ and the continuity properties in \eqref{eqn:outerapprox1} and \eqref{eqn:outerapprox2}, there exists $\rho \in (0,\min\{\epsilon/5,\epsilon/(5\kappa)]$ such that for all $(x',y',u',\alpha',\lambda') \in \ball(\bar w, \rho)$, 
\begin{align*}
  &\big|g^\nu(x',y') - g^\nu(\bar x, \bar y)\big| \leq \epsilon/5, ~~~~~~\big|g^\nu(x',z') - g^\nu(\bar x, z')\big| \leq \epsilon/5 ~~~\forall z' \in Y^\nu\\ 
  & \bar\lambda^\nu\Big|\dist\big(H^\nu(x,z'),D\big) - \dist\big(H^\nu(\bar x, z'),D\big)  \Big| \leq \epsilon/5~~~\forall z' \in Y^\nu.
\end{align*}
Since $|\alpha' - \bar \alpha| \leq \epsilon/5$ and $\kappa |\lambda' - \bar\lambda | \leq \epsilon/5$ when $(x',y',u',\alpha',\lambda') \in \ball(\bar w, \rho)$ and \eqref{eqn:psizexpan} holds, we find that for $w=(x,y,u,\alpha,\lambda)\in C \cap \ball(\bar w, \rho)$ and $z\in Y^\nu$ one has $|\psi_z(w) - \psi_z(\bar w)| \leq \epsilon$. 

For (v), maximizing $\psi_z(w)$ over $z\in Y^\nu$ is equivalent to minimizing $g^\nu(x,z) + \lambda \dist(H^\nu(x,z),D)$ over $z\in Y^\nu$, where $w = (x,y,u,\alpha,\lambda)$. Since $D$ is nonempty, $Y^\nu$ is nonempty and compact, and $g^\nu$ and $H^\nu$ are continuous, it follows via \cite[Theorem 4.9]{primer} that the minimum is attained.\eop

The assumptions of the proposition are mild. The conditions \eqref{eqn:outerapprox1} and \eqref{eqn:outerapprox2} simply amount to a uniformity across $z\in Y^\nu$ of the continuity properties for $g^\nu(\cdot, z)$ and $H^\nu(\cdot, z)$. The proposition applies even for $Y^\nu = Y$ but would then require that $Y$ is compact. 

{\blue 

\begin{example}{\rm (numerical results).}\label{ex:numerical} As a demonstration of the outer approximation algorithm, we consider three cases of {\rm (P)}$^\nu$ with lower-level problems in the form of integer linear programs and upper-level problems as linear programs. Let
\begin{align*}
f^\nu(x,y) & = \sum_{j=1}^n \bar c_j x_j  + \sum_{i=1}^m \bar a_i y_i, ~\mbox{ with } \bar c_j = -0.5(j-1),~  \bar a_i = -2(m-i+1)\\ 
g^\nu(x,y) & = \sum_{i=1}^m a_i y_i, ~\mbox{ with } a_i  = -0.5(i-1)\\
X & = \big\{x\in [0,1]^n~\big|~ A_0 x \leq b_0\big\}, ~Y = Y^\nu =\{0,1\}^m\\  
H^\nu(x,y) & = Ax + By - b,  ~D = (-\infty, 0]^m.
\end{align*}  
Table \ref{tab:num} summarizes numerical results for the outer approximation algorithm on three randomly generated cases with $n = m = q = 20$ and parameters $\sigma^\nu = \theta^\nu = \bar \lambda^\nu=10$ and $\tau^\nu=0$. A laptop with 1.7 GHz Intel processor and 16GB of RAM produces the computing times in the third column.     
\end{example}
\noindent {\bf  Detail}. The problem data is generated as follows. Each element in the $q\times n$-matrix $A$ is independently sampled from a uniform distribution on $[0,5]$, each element in the $q\times m$-matrix $B$ is independently sampled from a uniform distribution on $[0,10]$, and each element of the $r\times n$-matrix $A_0$ is independently sampled from a uniform distribution on $[0, 2]$; $r = 20$ in all cases. The vector $b$ has 30 as all its components and $b_0$ has 5 as its components. Adopting $\|\cdot\|_1$ as the norm and the approach described around \eqref{eqn:onenormexam} and \eqref{eqn:reformconstrs}, we obtain mixed-integer linear program formulations for the problems in both Step 1 and Step 2 after setting $\lambda = \bar\lambda^\nu$. We implement these formulations and the overall outer approximation algorithm in Python 3.12.1 using the Pyomo package \cite{Hartetal.17} and the Gurobi solver (version 12.0.0) with default options except MIPGap = 0.0001 and MIPGapAbs = 0. (Thus, the steps are accomplished with the tolerances $\delta^k$ and $\epsilon^k$ practically zero.) 

While Step 2 leads to a lower bound on the minimum value of (P)$^\nu$, Step 1 allows us to compute upper bounds after some straightforward calculations. The absolute value of the difference between these bounds divided by the lower bound is the {\em gap} reported in Table \ref{tab:num} after also multiplying by 100. Columns 4-8 in the table show how the gap decreases as the iterations of the outer approximation algorithm progress. {\bblue The occasional quick drop in gap appears to be caused by the finiteness of $Y^\nu$ and the linear structure in these instances.} We terminate the algorithm when the gap is no larger than 0.01\%; see the second column for the number of iterations.  The computational times vary across the instances, but this is expected since the algorithm is solving mixed-integer linear programs in both Step 1 and Step 2. A relaxed tolerance would produce shorter times (i.e., setting MIPGap $>$ 0.0001 and MIPGapAbs $>$ 0), and the same is likely to occur if the set $Y^\nu_k$ had been warm-started with a moderate set of points or if the Gurobi solver had been warm-started.\eop

}

\begin{table}[ht]

\begin{center}
\begin{tabular}{c|c|c|r|r|r|r|r|r}
         & iter. to    & time (sec.) to   & \multicolumn{5}{c}{~~~gap (\%) after $k$ iterations}\\
case & 0.01\% gap & 0.01\% gap & $k = 5$ & $k= 10$ & $k= 15$ & $k= 20$ & $k=25$\\ \hline
A    & 26   & 585  &  62.4   & 27.5    & 24.9    & 14.2    & 2.7 \\
B    & 17   & 80   &  56.8   & 52.4    & 13.7    & -       & - \\
C    & 27   & 630  &  51.8   & 49.3    & 44.7    & 28.1    & 25.2 
\end{tabular}\caption{Numerical results for the outer approximation algorithm across three randomly generated linear problem instances in Example \ref{ex:numerical}.}\label{tab:num}
\end{center}
\end{table}

{\bblue

\begin{example}{\rm (numerical results for quadratic problems).}\label{ex:numerical2} We also consider 27 mixed-integer quadratic instances of {\rm (P)}$^\nu$ by augmenting the construction in Example \ref{ex:numerical} with quadratic terms and adjusting the sets $X$ and $Y$. Let
\begin{align*}
f^\nu(x,y) & = \sum_{j=1}^n \bar c_j x_j  + \sum_{i=1}^m \bar a_i y_i + \sum_{j=1}^n x_j^2 + \sum_{i=1}^m y_i^2\\ 
g^\nu(x,y) & = \sum_{i=1}^m a_i y_i  + \sum_{j=1}^n(x_j-1)^2 + \sum_{i=1}^m(y_i-1)^2\\
X & = \big\{x\in X_0~\big|~ A_0 x \leq b_0\big\}, ~~Y = Y^\nu,
\end{align*}  
where $X_0$ is either $[0, 1]^n$ or $[0, 1]^{n/2} \times \{0,1\}^{n/2}$ or $\{0, 1\}^n$, with these cases referred to as ``real, ''r-b,'' and ``bin'' in Table \ref{tab:num2}, which reports optimality gaps for the outer approximation algorithm. Likewise, $Y$ is either $[0, 1]^m$ or $[0, 1]^{m/2} \times \{0,1\}^{m/2}$ or $\{0, 1\}^m$. With the three options $\bar \lambda^\nu=1$, $\bar \lambda^\nu=2$, and $\bar \lambda^\nu=10$, we obtain the 27 cases in Table \ref{tab:num2}. (The other components and parameters remain as in Example \ref{ex:numerical}.) 
\end{example}
\noindent {\bf  Detail}. Reformulations described around \eqref{eqn:onenormexam} and \eqref{eqn:reformconstrs} lead to mixed-integer (nonconvex) quadratic programs in Step 2 of the outer approximation algorithm for the given instances. Thus, we ease the computational burden by setting the MIPGap option for the Gurobi solver to 90\% of the current gap in the outer approximation algorithm. (Step 2 still provides lower bounds on the minimum value of (P)$^\nu$ by using Gurobi-provided lower bounds.) If Step 2 reproduces $x^k$ from the previous iteration or the computing time exceeds 900 seconds at Step 3, then the algorithm executes one final iteration with MIPGap option set to 3\% or 45\% of the current gap, whichever value is smaller. The final iteration is skipped if the current gap is already less than 1\%. The final iteration is terminated after 1800 seconds if the MIPGap tolerance is not reached. (The other implementation details remain as in Example \ref{ex:numerical}.) Table \ref{tab:num2} reports the resulting gaps (in percent) across the 27 instances. (The prevalence of gaps slightly below 3\% in the table is caused by the choice of \mbox{MIPGap} during the final iteration.) While the outer approximation algorithm achieves reasonable gaps when $y$ is binary, or $y$ is mixed and $x$ is binary, the algorithm suffers from weak lower bounds in other cases. Gaps tend to be smaller when $\bar \lambda^\nu$ is smaller because this parameter acts like a big-M in the formulation of Step 2. This motivates implementations that gradually increase this parameter as suggested by Assumption \ref{ass:basic}. We observe that when $y$ is real, then the semi-infinite constraints $\psi_z(w) \leq 0$, $z\in Y^\nu$, can be reformulated in the present instances using duality. This leads to a formulation of (P)$^\nu$ that might be solvable directly without outer approximations. Naturally, there are numerous algorithmic approaches for solving (P)$^\nu$ that might improve on the preliminary results reported here.\eop

}

\begin{table}[ht]
\begin{center}
\begin{tabular}{c||r|r|r||r|r|r||r|r|r}
 &  \multicolumn{3}{c||}{ $\bar\lambda^\nu = 1$}          &  \multicolumn{3}{c||}{ $\bar\lambda^\nu = 2$}  &  \multicolumn{3}{c}{ $\bar\lambda^\nu = 10$}\\
\hline
type        & bin $y$  & r-b $y$ & real $y$  & bin $y$  & r-b $y$ & real $y$ & bin $y$  & r-b $y$ & real $y$\\
\hline
bin  $x$      & 2.9   & 2.9   &  25.0        & 2.9  &  2.9   & 27.4     & 0.0  &  8.0    & 40.5\\ 
r-b $x$  & 2.9  & 8.6   &  39.0        & 0.8  & 24.3   & 50.0     & 1.5  & 51.1    & 68.4\\
real   $x$    & 1.2   & 13.8  &  49.1        & 0.5  & 27.8   & 46.6     & 0.6  & 50.3    & 63.0
\end{tabular}
\caption{{\bblue Percent optimality gaps when applying the outer approximation algorithm to 27 nonlinear problem instances in Example \ref{ex:numerical2} involving binary (bin), real and binary (r-b), and real (real) variables.}}\label{tab:num2}
\end{center}
\end{table}

\section{Intermediate Results and Proofs}\label{sec:proofs}

In this section, we use the notation: For $x\in\reals^n$ and $\lambda\in [0,\infty)$, let 
\begin{align*}
\mu(x,\lambda) &= \inf_{y\in Y} g(x,y) + \lambda \dist\big(H(x,y), D\big)\\
\mu^\nu(x,\lambda) &= \inf_{y\in Y^\nu} g^\nu(x,y) + \lambda \dist\big(H^\nu(x,y), D\big).
\end{align*}

\begin{lemma}{\rm (lower approximation of penalty-based value functions).}\label{lemma:lbmu} Suppose that Assumption \ref{ass:basic}(a-d) holds. If $x^\nu\in X\to x$ and $\lambda^\nu \in [0,\infty) \to \lambda$, then $\nlimsup \mu^\nu(x^\nu,\lambda^\nu) \leq \mu(x,\lambda)$.
\end{lemma}
\noindent {\bf  Proof}. Let $\gamma \in (0,\infty)$. We consider two cases. (i) If $\mu(x,\lambda) \in \reals$, then there exists $y\in Y$ such that 
\begin{equation}\label{eqn:part1limsupres}
g(x, y) + \lambda \dist\big(H(x, y), D\big) \leq \mu(x,\lambda) + \gamma. 
\end{equation}
Since $Y^\nu\sto Y$, there exists a sequence $\{y^\nu\in Y^\nu, \nu\in\nats\}$ converging to $y$. The continuity of $g$ relative to $X\times Y$ and the vanishing errors in Assumption \ref{ass:basic}(c) yield 
\begin{equation}\label{eqn:part2limsupres}
\nlimsup g^\nu(x^\nu,y^\nu) \leq \nlimsup \big(g^\nu(x^\nu,y^\nu) - g(x^\nu,y^\nu)\big) + \nlimsup g(x^\nu,y^\nu) = g(x,y). 
\end{equation}
Similarly, by Assumption \ref{ass:basic}(d) and continuity of $\dist(H(\cdot,\cdot), D)$ relative to $X\times Y$, one has 
\begin{equation}\label{eqn:part3limsupres}
\nlimsup \dist\big(H^\nu(x^\nu,y^\nu), D\big) \leq \dist\big(H(x,y), D\big). 
\end{equation}
Combining the inequalities in \eqref{eqn:part1limsupres}-\eqref{eqn:part3limsupres} while also incorporating $\lambda^\nu$ and $\lambda$, we obtain 
\begin{align*}
\nlimsup \mu^\nu(x^\nu,\lambda^\nu) & \leq \nlimsup \Big(g^\nu(x^\nu,y^\nu) + \lambda^\nu \dist\big(H^\nu(x^\nu,y^\nu), D\big)\Big)\\
& \leq \nlimsup g^\nu(x^\nu,y^\nu) + \nlimsup \Big(\lambda^\nu \dist\big(H^\nu(x^\nu,y^\nu), D\big)\Big)\\
& \leq g(x,y) + \lambda \dist\big(H(x,y), D\big) \leq \mu(x,\lambda) + \gamma. 
\end{align*}

(ii) If $\mu(x,\lambda) = - \infty$, then there exists $y\in Y$ such that $g(x, y) + \lambda \dist\big(H(x, y), D\big)$ $\leq$ $-1/\gamma$. Again constructing $\{y^\nu\in Y^\nu, \nu\in\nats\}$ and leveraging \eqref{eqn:part2limsupres}-\eqref{eqn:part3limsupres}, we obtain 
\[
\nlimsup \mu^\nu(x^\nu,\lambda^\nu) \leq g(x,y) + \lambda \dist\big(H(x,y), D\big)\leq -1/\gamma. 
\]
Since $\gamma$ is arbitrary, the conclusion holds in both cases.\eop

\begin{lemma}{\rm (vanishing error for penalty-based value function).}\label{lemma:convmu}
For nonempty compact set $Y\subset\reals^m$, nonempty closed set $D\subset\reals^q$, $g:\reals^n\times\reals^m\to \reals$, and $H:\reals^n\times\reals^m\to \reals^q$, suppose that $x\in \reals^n$ is fixed and $g(x,\cdot)$ and $H(x,\cdot)$ are continuous relative to $Y$. Then, the following hold: 
\begin{enumerate}[(a)]
\item If $\lambda^\nu \to \infty$, then $\nlim \mu(x,\lambda^\nu) = \inf_{y\in Y} \{g(x,y)~|~H(x,y) \in D\}$.

\item If $\tau \in (0,\infty)$ and $\inf_{y\in Y} \{g(x,y)~|~H(x,y) \in D\}<\infty$, then there exist $\bar y\in Y$ and $\lambda\in [0,\infty)$ such that $H(x,\bar y) \in D$ and $g(x,\bar y) \leq \mu(x,\lambda) + \tau$.
\end{enumerate}  
\end{lemma}
\noindent {\bf Proof}. Let $\psi,\psi^\nu:\reals^m\to (-\infty,\infty]$ be given by 
\[
\psi(y) = \iota_Y(y) + g(x,y) + \iota_D\big( H(x,y) \big); ~~~~\psi^\nu(y) = \iota_Y(y) + g(x,y) + \lambda^\nu \dist\big(H(x,y), D\big).
\]
For part (a), we consider two cases. (i) Suppose that $H(x,y) \not \in D$ for all $y\in Y$. For the sake of contradiction, suppose that there exist $\gamma < \infty$ and a subsequence $N'\subset\nats$ such that $\mu(x,\lambda^\nu) \leq \gamma$ for all $\nu\in N'$. 
For each $\nu$, there exists $y^\nu \in \nargmin_{y\in \reals^m} \psi^\nu(y)$ because $Y$ is compact and nonempty and $g(x,\cdot)$ and $H(x,\cdot)$ are continuous relative to $Y$; see, e.g., \cite[Theorem 4.9]{primer}. By the compactness of $Y$, there exist a point $\bar y\in Y$ and a subsequence $N\subset N'$ such that $y^\nu \Nto \bar y$. Since $H(x,\cdot)$ is continuous relative to $Y$, this implies that $\dist(H(x,y^\nu),D) \Nto \dist(H(x,\bar y), D)$. The set $D$ is closed and therefore $\dist(H(x,\bar y), D) > 0$ and $\{\dist(H(x,y^\nu),D), \nu\in N\}$ is bounded away from zero. Moreover, $g(x,\cdot)$ is continuous relative to $Y$, which implies that $\{g(x,y^\nu), \nu\in N\}$ is bounded from below. Combining these facts, we conclude that 
$\mu(x,\lambda^\nu) = g(x,y^\nu) + \lambda^\nu \dist(H(x,y^\nu),D) \Nto \infty$, which is a contradiction. 

(ii) Suppose that $H(x,\bar y) \in D$ for some $\bar y\in Y$. Then, $\inf \psi < \infty$. Trivially, we have $\psi(y) > - \infty$ for all $y\in \reals^m$. Since $Y$ is compact, $\{\psi^\nu, \nu\in \nats\}$ is tight. If $\psi^\nu\eto \psi$, then $\inf \psi^\nu \to \inf \psi$ by \cite[Theorem 5.5(d)]{primer}. Thus it only remains to verify such epi-convergence. 

First, we verify the liminf-condition \eqref{eqn:liminf}. Suppose that $y^\nu\to y$. Without loss of generality, we assume that $y^\nu \in Y$ because otherwise $\psi^\nu(y^\nu) = \infty$. Since $Y$ is closed, this implies that $y\in Y$. If $H(x,y) \in D$, then it follows by continuity of $g(x,\cdot)$ relative to $Y$ that 
\[
\nliminf \psi^\nu(y^\nu) \geq \nliminf g(x,y^\nu) + \nliminf \lambda^\nu \dist\big(H(x,y^\nu), D\big) \geq g(x,y) = \psi(y).  
\]
If $H(x,y) \not\in D$, then there exists $\epsilon > 0$ such that $\dist(H(x,y),D) = \epsilon$ because $D$ is closed. The continuity of $H(x,\cdot)$ relative to $Y$ implies that $\dist(H(x,y^\nu), D) \geq \epsilon/2$ for all $\nu$ sufficiently large. Thus, $\psi^\nu(y^\nu) \to \infty$ and again 
$\nliminf \psi^\nu(y^\nu) \geq \psi(y)$.

Second, we verify the limsup-condition \eqref{eqn:limsup}. Let $y \in \reals^m$. Without loss of generality, we assume that $y\in Y$ and $H(x,y)\in D$ because otherwise $\psi(y) = \infty$. We construct $y^\nu = y$. Then, 
\[
\nlimsup \psi^\nu(y^\nu) \leq g(x,y) + \nlimsup \lambda^\nu \dist\big(H(x,y), D\big) = g(x,y) = \psi(y).  
\]  
We conclude that $\psi^\nu\eto \psi$.

For part (b), we observe that $\inf \psi>-\infty$ because $Y$ is compact, $D$ is closed, and $g(x,\cdot)$ and $H(x,\cdot)$ are continuous relative to $Y$. Moreover, $\inf \psi < \infty$ be assumption. Thus, there exists $\bar y\in Y$ such that $H(x,\bar y) \in D$ and $g(x,\bar y) \leq \inf \psi + \tau/2$. By part (a), there is $\lambda \in [0,\infty)$ such that $\mu(x,\lambda) \geq \inf \psi - \tau/2$. Consequently, $g(x,\bar y) \leq \inf \psi + \tau/2 \leq \mu(x,\lambda) + \tau$ as claimed.\eop

Lemma \ref{lemma:convmu}(a) is false if the compactness of $Y$ is removed. For a counterexample with a convex objective function and a convex feasible set, let $Y = \reals$, $g(x,y) = -y$, $D= (-\infty, 0]$, and $H(x,y) = \min\{y+1, \exp(-y)\}$. Then, $\mu(x,\lambda) = -\infty$ for any $\lambda \in [0,\infty)$ but $\inf_{y\in Y} \{g(x,y)~|~H(x,y) \in D\} = 1$.

For the {\em remainder of the paper}, we represent (P)$^\nu$ using the extended real-valued function $\phi^\nu:\reals^n\times\reals^m\times \reals^q \times \reals \times \reals \to (-\infty, \infty]$ defined by 
\begin{align*}
\phi^\nu(x,y,u,\alpha,\lambda) & = \iota_X(x) + \iota_Y(y) + f^\nu(x,y) + \sigma^\nu \|u\| - \theta^\nu \alpha + \iota_{D}\big(H^\nu(x,y) + u\big)\\
& ~~~+ \iota_{(-\infty, \tau^\nu]}\big( g^\nu(x,y) + \alpha - \mu^\nu(x,\lambda) \big) + \iota_{(-\infty,0]}(\alpha) + \iota_{[0,\bar\lambda^\nu]}(\lambda).
\end{align*}
It is apparent $\dom \phi^\nu$ coincides with the feasible set of (P)$^\nu$ because $g^\nu(x,y) + \alpha - \mu^\nu(x,\lambda) \leq \tau^\nu$ if and only if \eqref{eqn:PnuConstr_mu} holds. On that feasible set, $\phi^\nu$ matches the objective function of (P)$^\nu$. Thus, the set of $\epsilon$-optimal solutions of (P)$^\nu$ coincides with $\epsilon$-$\nargmin \phi^\nu$ and $\mathfrak{m}^\nu= \inf \phi^\nu$. 

For any $\beta \in [0,\infty)$, we also define $\phi_\beta:\reals^n\times\reals^m\times \reals^q \times \reals \times \reals \to (-\infty, \infty]$ by setting 
\begin{align*}
\phi_\beta(x,y,u,\alpha,\lambda) & = \iota_X(x) + \iota_Y(y) + f(x,y) + \iota_{D}\big(H(x,y) + u\big) \\
&~~~ + \iota_{(-\infty, \beta]}\big( g(x,y) + \alpha - \mu(x,\lambda) \big) + \iota_{\{0\}^q}(u) + \iota_{\{0\}}(\alpha) + \iota_{[0,\infty)}(\lambda).
\end{align*}
If $\beta = \tau$, then minimizing $\phi_\beta$ is a restriction of (P) because of the reformulation \eqref{eqn:level-reform} and the fact that $\mu(x,\lambda) \leq \inf_{y\in Y} \{g(x,y) \, |\, H(x,y) \in D\}$ for any $x\in \reals^n$ and $\lambda \in [0,\infty)$. 

The next two lemmas establish relationships between $\phi^\nu$ and $\phi_\beta$. When combined, they yield the epi-convergence $\phi^\nu \eto \phi_\beta$ under natural assumptions.

\begin{lemma}\label{lemma:liminf}{\rm (liminf property).} Suppose that Assumption \ref{ass:basic}(a-d) holds and the functions $f,f^\nu:\reals^n\times\reals^m\to (-\infty,\infty]$ satisfy $\nliminf f^\nu(x^\nu,y^\nu) \geq f(x,y)$  when $x^\nu \in X\to x$ and $y^\nu \in Y\to y$.
If $\tau^\nu \in [0,\infty) \to \beta$ and $\sigma^\nu, \theta^\nu \to \infty$, then 
\begin{equation*}
\nliminf \phi^\nu(x^\nu,y^\nu,u^\nu,\alpha^\nu,\lambda^\nu) \geq \phi_\beta(x,y,u,\alpha,\lambda) \mbox{ for } (x^\nu,y^\nu,u^\nu,\alpha^\nu,\lambda^\nu)\to (x,y,u,\alpha,\lambda).
\end{equation*}
\end{lemma}
\noindent {\bf Proof.} Let $(x^\nu,y^\nu,u^\nu,\alpha^\nu,\lambda^\nu)\to (x,y,u,\alpha,\lambda)$. Without loss of generality, we assume that $x^\nu\in X$, $y^\nu\in Y$, $\alpha^\nu\leq 0$, $\lambda^\nu \in [0, \bar \lambda^\nu]$ for all $\nu$ because other values of $(x^\nu,y^\nu,\alpha^\nu,\lambda^\nu)$ would cause $\phi^\nu(x^\nu,y^\nu,u^\nu,\alpha^\nu,\lambda^\nu) = \infty$ regardless of $u^\nu$, which is unproblematic in the inequality we seek to establish. After invoking closedness, we proceed under the assumption that $x\in X, y\in Y, \alpha \leq 0$, and $\lambda \in [0,\infty)$. We consider five cases: 

(i) Suppose that $u\neq 0$ so that $\phi_\beta(x,y,u,\alpha,\lambda) = \infty$. Since $u^\nu\to u$ and $\sigma^\nu \to \infty$, we find that $\sigma^\nu \|u^\nu\| \to \infty$. Thus, $\phi^\nu(x^\nu,y^\nu,u^\nu,\alpha^\nu,\lambda^\nu) \to \infty$ because $\nliminf f^\nu(x^\nu,y^\nu) \geq f(x,y) > - \infty$. 

(ii) Suppose that $\alpha < 0$ so that $\phi_\beta(x,y,u,\alpha,\lambda) = \infty$. Since $\alpha^\nu\to \alpha$ and $\theta^\nu \to \infty$, we find that $-\theta^\nu \alpha^\nu \to \infty$. Thus, again, $\phi^\nu(x^\nu,y^\nu,u^\nu,\alpha^\nu,\lambda^\nu) \to \infty$. 

(iii) Suppose that $u = 0$, $\alpha = 0$, $H(x,y)\in D$, and $g(x,y) - \mu(x,\lambda) \leq \beta$. Then, one has 
\begin{align*}
&\nliminf \, \phi^\nu(x^\nu,y^\nu,u^\nu,\alpha^\nu,\lambda^\nu)  \geq \nliminf \iota_X(x^\nu) + \nliminf \iota_Y(y^\nu) + \nliminf f^\nu(x^\nu,y^\nu)\\
& + \nliminf \sigma^\nu \|u^\nu\| + \nliminf (-\theta^\nu \alpha^\nu) + \nliminf \iota_{D}\big(H^\nu(x^\nu,y^\nu) + u^\nu\big) + \nliminf \iota_{[0,\bar\lambda^\nu]}(\lambda^\nu)\\
& + \nliminf \iota_{(-\infty, \tau^\nu]}\big( g^\nu(x^\nu,y^\nu) + \alpha^\nu - \mu^\nu(x^\nu,\lambda^\nu)\big) + \nliminf \iota_{(-\infty,0]}(\alpha^\nu)\\
& \geq f(x,y) = \phi_\beta(x,y,u,\alpha,\lambda). 
\end{align*}

(iv) Suppose that $u = 0$, $\alpha = 0$, and $H(x,y) \not\in D$. Then, $\phi_\beta(x,y,u,\alpha,\lambda) = \infty$. Since $D$ is closed, there is $\epsilon > 0$ such that $\dist(H(x,y),D) = \epsilon$.  We recall that $u^\nu\to 0$, $\alpha^\nu\to 0$, and Assumption \ref{ass:basic}(d) holds. Thus, there is $\bar \nu$ such that $\|u^\nu\| \leq \epsilon/4$ and
\[
\dist\big(H(x^\nu,y^\nu), D\big) \geq 3\epsilon/4;  \dist\big( H^\nu(x^\nu,y^\nu), D\big) \geq \dist\big( H(x^\nu,y^\nu), D\big) - \epsilon/4 ~\forall \nu\geq \bar\nu.
\]
Consequently, for $\nu\geq \bar\nu$, 
\begin{align*}
&\dist\big( H^\nu(x^\nu,y^\nu) + u^\nu, D\big)   \geq \dist\big( H^\nu(x^\nu,y^\nu), D\big) - \|u^\nu\|\\
& \geq \dist\big( H(x^\nu,y^\nu), D\big) - \epsilon/4 -\epsilon/4 \geq 3\epsilon/4 - \epsilon/4 - \epsilon/4  = \epsilon/4 > 0.  
\end{align*}
It follows that $\phi^\nu(x^\nu,y^\nu,u^\nu,\alpha^\nu,\lambda^\nu) = \infty$ for all $\nu\geq \bar \nu$.

(v) Suppose that $u = 0$, $\alpha = 0$, and there is $\epsilon \in (0,\infty)$ such that $g(x,y) - \mu(x,\lambda) = \beta + \epsilon$. Then,  $\phi_\beta(x,y,u,\alpha,\lambda) = \infty$. By invoking Lemma \ref{lemma:lbmu}, Assumpton \ref{ass:basic}(c), and the facts that $\alpha^\nu\to 0$ and $\tau^\nu\to \beta$, there exists $\bar \nu$ such that for all $\nu\geq \bar\nu$: 
\begin{align*}
& \mu^\nu(x^\nu,\lambda^\nu) \leq \mu(x,\lambda) + \epsilon/6; ~~~~ g(x^\nu,y^\nu) \geq g(x,y) - \epsilon/6;  ~~~~ \alpha^\nu \geq -\epsilon/6\\ 
& g^\nu(x^\nu,y^\nu) \geq g(x^\nu,y^\nu) - \epsilon/6; ~~~~ \tau^\nu \leq \beta + \epsilon/6.
\end{align*}
Thus, for $\nu\geq \bar\nu$, 
\begin{align*}
g^\nu(x^\nu,y^\nu) + \alpha^\nu - \mu^\nu(x^\nu,\lambda^\nu) & \geq g(x^\nu,y^\nu) - \epsilon/6 - \epsilon/6 - \mu(x,\lambda) - \epsilon/6\\
& \geq g(x,y) - \epsilon/6 - \epsilon/6 - \epsilon/6 - \mu(x,\lambda) - \epsilon/6\\
& = \beta + \epsilon - 4\epsilon/6 \geq \tau^\nu + \epsilon/6.  
\end{align*}
We conclude that $\phi^\nu(x^\nu,y^\nu,u^\nu,\alpha^\nu,\lambda^\nu) = \infty$ for all $\nu\geq \bar \nu$.\eop

\begin{lemma}\label{lemma:limsup}{\rm (limsup property).}
If Assumption \ref{ass:basic} holds, $\beta \in [0,\infty)$, $\{\theta^\nu, \nu\in\nats\}$ is bounded away from zero, and $\theta^\nu|\tau^\nu-\beta|\to 0$, then, for any $(x,y,u,\alpha,$ $\lambda)$ $\in$ $\reals^n\times\reals^m\times\reals^q\times\reals\times\reals$, there exists a sequence $(x^\nu,y^\nu,u^\nu,\alpha^\nu,\lambda^\nu)\to (x,y,u,\alpha,\lambda)$ with
\begin{equation*}
\nlimsup \phi^\nu(x^\nu,y^\nu,u^\nu,\alpha^\nu,\lambda^\nu) \leq \phi_\beta(x,y,u,\alpha,\lambda).
\end{equation*}

\end{lemma}
\noindent {\bf Proof.} Let $x \in \reals^n$, $y\in\reals^m$, $u\in \reals^q$, $\alpha\in\reals$, and $\lambda\in \reals$. Without loss of generality, we assume that $x\in X$, $y\in Y$, $u = 0$, $\alpha = 0$, $H(x,y) \in D$, $g(x,y) - \mu(x,\lambda) \leq \beta$, and $\lambda \geq 0$ because otherwise $\phi_\beta(x,y,u,\alpha,\lambda) = \infty$ and the conclusion holds trivially with $(x^\nu,y^\nu,u^\nu,\alpha^\nu,\lambda^\nu) = (x,y,u,\alpha,\lambda)$. 

We construct $x^\nu = x$, $y^\nu = y$, $\lambda^\nu = \min\{\bar\lambda^\nu, \lambda\}$, and 
\begin{align*}
  u^\nu & \in \nargmin_{\bar u\in\reals^q} \big\{\|\bar u\|~\big|~H^\nu(x,y) + \bar u \in D\big\}\\
  \alpha^\nu & = \min\big\{0, g(x,y) - \mu(x,\lambda) - g^\nu(x,y) + \mu^\nu(x,\lambda^\nu) - \beta + \tau^\nu\big\}.
\end{align*}
Since $D$ is nonempty and closed, $u^\nu$ is well defined. The construction implies that $H^\nu(x^\nu,y^\nu) + u^\nu \in D$. The condition $g(x,y) - \mu(x,\lambda) \leq \beta$ ensures that $\mu(x,\lambda)\in \reals$. Thus, $\alpha^\nu$ is well-defined but with value $-\infty$ if $\mu^\nu(x,\lambda^\nu) = -\infty$. Trivially, $\alpha^\nu \leq 0$ and $\lambda^\nu \in [0, \bar\lambda^\nu]$. 

Next, we verify that $\sigma^\nu \|u^\nu\| \to 0$. Invoking Assumption \ref{ass:basic}(d), we see that
\[
\|u^\nu\| = \dist\big(H^\nu(x,y), D\big) \leq \dist\big(H(x,y), D\big) + \eta^\nu \leq \eta^\nu. 
\]
Thus, $\|u^\nu\| \to 0$ and $\sigma^\nu \|u^\nu\| \to 0$ follow via Assumption \ref{ass:basic}(d) and \ref{ass:basic}(f), respectively.  

To establish $-\theta^\nu \alpha^\nu \to 0$, we invoke Assumption \ref{ass:basic}(c) to obtain  
\begin{align}
- \alpha^\nu & = \max\big\{0, \mu(x,\lambda)  - \mu^\nu(x,\lambda^\nu) + g^\nu(x,y) - g(x,y) + \beta - \tau^\nu\big\}\nonumber\\
          & \leq \max\big\{0, \mu(x,\lambda)  - \mu^\nu(x,\lambda^\nu) + \delta^\nu + |\beta - \tau^\nu|\big\}.\label{eqn:u2bd}  
\end{align}
The remaining difference on the right-hand side in this inequality is bounded: 
\begin{align}
\mu(x,\lambda)  - \mu^\nu(x,\lambda^\nu) &\leq \inf_{z\in Y} \Big\{ g(x,z) + \lambda \dist\big( H(x,z), D\big) \Big\}\label{eqn:mubdpart2}\\
& - \inf_{z\in Y} \Big\{ g^\nu(x,z) + \lambda^\nu \dist\big( H^\nu(x,z), D\big) \Big\},\nonumber
\end{align}
where we leverage the fact that $Y^\nu \subset Y$. By Assumption \ref{ass:basic}(c,d),  
\[
g(x,z) + \lambda \dist\big(H(x,z), D\big) - g^\nu(x,z) - \lambda\dist\big(H^\nu(x,z), D \big) \leq \delta^\nu + \lambda \eta^\nu \leq \delta^\nu + \bar\lambda^\nu \eta^\nu
\]
regardless of $z\in Y$. We can use this fact in conjunction with \eqref{eqn:mubdpart2} because $\lambda^\nu = \lambda$ for sufficiently large $\nu$ by construction and Assumption \ref{ass:basic}(f). Thus, for sufficiently large $\nu$, $\mu(x,\lambda)  - \mu^\nu(x,\lambda^\nu) \leq \delta^\nu + \bar\lambda^\nu \eta^\nu$. We combine this fact with \eqref{eqn:u2bd} to establish that $- \alpha^\nu \leq 2\delta^\nu + \bar\lambda^\nu \eta^\nu + |\beta - \tau^\nu|$ for sufficiently large $\nu$. Since $-\alpha^\nu \geq 0$ by construction, it follows that $-\theta^\nu \alpha^\nu  \to 0$ by Assumption \ref{ass:basic}(f) and the assumed rate by which $\tau^\nu\to \beta$. 

Since $\mu^\nu(x,\lambda^\nu)$ is finite for sufficiently large $\nu$, the definition of $\alpha^\nu$ implies that   
\begin{align*}
&  g^\nu(x^\nu,y^\nu) + \alpha^\nu - \mu^\nu(x^\nu,\lambda^\nu)\\
& \leq g^\nu(x,y) + g(x,y) - \mu(x,\lambda) - g^\nu(x,y) + \mu^\nu(x,\lambda^\nu) - \beta + \tau^\nu - \mu^\nu(x,\lambda^\nu)\\ 
& \leq \beta - \beta + \tau^\nu = \tau^\nu
\end{align*}
for such $\nu$. Collecting these facts, we confirm that 
\begin{align*}
 & \nlimsup \phi^\nu(x^\nu,y^\nu,u^\nu,\alpha^\nu,\lambda^\nu) \leq \nlimsup ( - \theta^\nu \alpha^\nu) + \nlimsup f^\nu(x,y) + \nlimsup \sigma^\nu \|u^\nu\|\\
 & +\iota_X(x)  + \nlimsup \iota_{D}\big(H^\nu(x,y) + u^\nu\big) + \nlimsup \iota_{(-\infty, \tau^\nu]}\big( g^\nu(x,y) + \alpha^\nu - \mu^\nu(x,\lambda^\nu)\big)\\
 &  + \iota_Y(y) + \nlimsup \iota_{(-\infty,0]}(\alpha^\nu) + \nlimsup \iota_{[0,\bar\lambda^\nu]}(\lambda^\nu) \leq f(x,y) = \phi_\beta(x,y,u,\alpha,\lambda), 
\end{align*}
where Assumption \ref{ass:basic}(e) also enters. It is apparent that $(x^\nu,y^\nu,u^\nu,\alpha^\nu,\lambda^\nu)$ $\to$ $(x,y,u,$ $\alpha,\lambda)$; $\alpha^\nu \to 0$ because $\theta^\nu \alpha^\nu \to 0$ and $\{\theta^\nu, \nu\in\nats\}$ is bounded away from zero. Thus, the conclusion follows.\eop

The liminf-condition of Assumption \ref{ass:basic}(e) is actually superfluous in the lemma and can be omitted.

\medskip

\noindent {\bf Proof of Theorem \ref{thm:withoutcalmness}(a)}. We invoke Lemma \ref{lemma:limsup}, with $\beta = \tau'$. Thus, for all 
$(x,y,u,\alpha,\lambda)\in\reals^n\times\reals^m\times\reals^q\times\reals\times\reals$, there is a sequence $(x^\nu,y^\nu,u^\nu,\alpha^\nu,\lambda^\nu)\to (x,y,u,\alpha,\lambda)$ with 
\[
\nlimsup \phi^\nu(x^\nu,y^\nu,u^\nu,\alpha^\nu,\lambda^\nu) \leq \phi_{\tau'}(x,y,u,\alpha,\lambda).
\]
This fact allows us to follow the proof of \cite[Theorem 5.5(a)]{primer} and deduce that $\nlimsup(\inf \phi^\nu)$ $\leq$ $\inf \phi_{\tau'}$. Next, let $\gamma \in (0, \tau' - \tau]$. We consider three cases: 

(i) Suppose that $\mathfrak{m} \in \reals$. Then, there are $x \in X$ and $y\in \tau\mbox{-}\nargmin_{z\in Y} \{g(x,z)~|~$ $H(x,z) \in D\}$ such that $f(x,y) \leq \mathfrak{m} + \gamma$. By \eqref{eqn:level-reform}, $y\in Y$, $H(x,y) \in D$, and 
\[
g(x,y) \leq \inf_{z\in Y} \big\{g(x,z)~\big|~H(x,z) \in D\big\} + \tau.
\]
Lemma \ref{lemma:convmu}(a) ensures that there exists $\hat \lambda \in [0,\infty)$ such that 
\[
\mu(x,\hat\lambda) \geq \inf_{z\in Y} \big\{g(x,z)~\big|~H(x,z) \in D\big\} - \gamma. 
\]
Set $\alpha = 0$. Then, 
\[
g(x,y) + \alpha \leq \inf_{z\in Y} \big\{g(x,z)~\big|~H(x,z) \in D\big\} + \tau \leq \mu(x,\hat \lambda) + \gamma + \tau \leq \mu(x,\hat \lambda) + \tau'.
\]
Since $\mathfrak{m}^\nu = \inf \phi^\nu$ by construction, it follows that
\[
\nlimsup \mathfrak{m}^\nu = \nlimsup\big(\inf \phi^\nu\big) \leq \inf \phi_{\tau'} \leq \phi_{\tau'}(x,y,0,0,\hat \lambda) = f(x,y) \leq \mathfrak{m} + \gamma. 
\]
Since $\gamma$ can be made arbitrarily close to zero, the claim follows in this case. 

(ii) Suppose that $\mathfrak{m} = -\infty$. Again passing through \eqref{eqn:level-reform}, we realize that there exist $x \in X$ and $y\in Y$ such that $H(x,y) \in D$, $g(x,y) \leq \inf_{z\in Y} \{g(x,z)~|~H(x,z) \in D\} + \tau$ and $f(x,y) \leq -1/\gamma$. Mimicking the arguments of case (i), we find that $\nlimsup \mathfrak{m}^\nu \leq -1/\gamma$ and the claim holds in this case too because $\gamma$ can be chosen arbitrarily close to zero.

(iii) If $\mathfrak{m} = \infty$, then the claim holds trivially.\eop

\noindent {\bf Proof of Theorem \ref{thm:withoutcalmness}(b)}. Lemmas \ref{lemma:liminf} and \ref{lemma:limsup} hold with $\beta = \tau$. This implies that $\phi^\nu \eto \phi_\tau$; see \eqref{eqn:liminf} and \eqref{eqn:limsup}. The existence of a compact set $B$ implies that $\{\phi^\nu, \nu\in\nats\}$ is tight. We then bring in \cite[Theorem 5.5(d)]{primer} to conclude that $\inf \phi^\nu \to \inf \phi_\tau$ provided that $\inf \phi_\tau<\infty$ and $\phi_\tau$ takes nowhere the value $-\infty$. While the latter is easily ruled out, $\inf \phi_\tau = \infty$ is possible under the current assumptions. However, we can step through the proof of \cite[Theorem 5.5(d)]{primer} and verify that $\inf \phi^\nu \to \inf \phi_\tau$ holds even in this case. Consequently, 
$\mathfrak{m}^\nu = \inf \phi^\nu \to \inf \phi_\tau \geq \mathfrak{m}$, 
where the last inequality follows from the fact that minimizing $\phi_\tau$ is a restriction of (P) due to \eqref{eqn:level-reform} and the relationship $\mu(x,\lambda)$ $\leq$ $\inf_{z\in Y} \{f(x,z)~|~H(x,z) \in D\}$ for all $x\in\reals^n$ and $\lambda \in [0,\infty)$.\eop 

\noindent {\bf Proof of Theorem \ref{thm:withoutcalmness}(c)}. Lemmas \ref{lemma:liminf} and \ref{lemma:limsup} hold with $\beta = \tau$. This implies that $\phi^\nu \eto \phi_\tau$. The assumption in \eqref{eqn:extraass} about the existence of a point $(x,y)$ with certain properties ensures that $\inf \phi_\tau<\infty$. Since $f$ takes nowhere the value $-\infty$, $\phi_\tau$ is also greater than $-\infty$. We can then invoke \cite[Theorem 5.5(b)]{primer} to establish that 
$(\hat x, \hat y, \hat u, \hat \alpha, \hat \lambda) \in \epsilon\mbox{-}\nargmin \phi_\tau$, where $\epsilon = \sup_{\nu\in\nats} \epsilon^\nu$, which is finite by assumption. 
Thus, $\phi_\tau(\hat x, \hat y, \hat u, \hat \alpha, \hat \lambda)<\infty$, which implies that $\hat x \in X$, $\hat y \in Y$, $\hat u = 0$, $\hat \alpha = 0$, $H(\hat x, \hat y) \in D$, and 
\[
g(\hat x, \hat y) \leq \mu(\hat x, \hat \lambda) + \tau \leq \inf_{z\in Y} \big\{g(\hat x, z)~\big|~H(\hat x, z) \in D\big\} + \tau.
\] 
The conclusion follows by \eqref{eqn:level-reform}.

The claim about \eqref{eqn:extraass} holds because if (P) has a feasible point, then there exist $x\in X$ and $y\in Y$ such that $H(x,y)\in D$. Thus, $\inf_{z\in Y} \{g(x,z)~|~H(x,z)\in D\} < \infty$ and we invoke Lemma \ref{lemma:convmu}(b).\eop

\noindent {\bf Proof of Proposition \ref{prop:char_calm}}. If $Y$ is empty, then the claim holds trivially. Thus, we proceed under the assumption that $Y$ is nonempty. First, suppose that $V$ is calm at $x$ with penalty threshold $\lambda$. Let $y\in Y$ and $\epsilon >0$. Since $D$ is nonempty, there exists $\bar u\in D$ such that 
\[
\big\|H(x,y) - \bar u\big\| \leq \dist\big(H(x,y), D\big) + \epsilon. 
\]
With $u = \bar u - H(x,y)$, one has $H(x,y) + u \in D$ and $\|u\| \leq \dist(H(x,y), D) + \epsilon$. Calmness yields
\begin{align*}
V(x,0) & \leq V(x,u) + \lambda \|u\|  \leq V(x,u) + \lambda\dist\big(H(x,y), D\big) + \lambda\epsilon\\
& = \inf_{z\in Y} \Big\{ g(x,z)~\Big|~ H(x,z) + u \in D \Big\} + \lambda\dist\big(H(x,y), D\big) + \lambda\epsilon\\
 & \leq g(x,y) + \lambda\dist\big(H(x,y), D\big) + \lambda\epsilon. 
\end{align*}
Since $\epsilon$ and $y$ are arbitrary, we have established that $V(x,0) \leq \mu(x,\lambda)$. The reverse inequality $V(x,0) \geq \mu(x,\lambda)$ holds trivially.  

Second, suppose that $\mu(x,\lambda) = V(x,0)$. Let $y\in Y$ and $u\in\reals^q$ satisfy $H(x,y) + u \in D$. Then, 
\begin{align*}
  V(x,0) = \mu(x,\lambda) & \leq g(x,y) + \lambda \dist\big(H(x,y), D \big)\\
& \leq g(x,y) + \lambda \dist\big(H(x,y)+u, D \big) + \lambda \|u\| = g(x,y) + \lambda\|u\|. 
\end{align*} 
Since this holds for arbitrary $y\in Y$ satisfying $H(x,y) + u \in D$, we also have $V(x,0) \leq V(x,u) + \lambda\|u\|$ and the conclusion follows.\eop

\begin{lemma}\label{lemma:phi}{\rm (properties of $\phi_\tau$).} Suppose there exists an optimal solution $(x^\star,y^\star)$ of {\rm (P)} such that the value function $V$ is calm at $x^\star$. Then, the following hold:
\begin{enumerate}[(a)]
\item $\inf \phi_\tau = \mathfrak{m}\in\reals$.

\item $(\bar x, \bar y, \bar u, \bar \alpha, \bar \lambda) \in \nargmin \phi_\tau$ implies that $(\bar x, \bar y)$ is an optimal solution of {\rm (P)}.

\item If $\lambda^\star$ is a penalty threshold for the value function $V$ at $x^\star$, then $(x^\star, y^\star, 0, 0, \lambda) \in \nargmin \phi_\tau$ for any $\lambda \in [\lambda^\star, \infty)$. 
\end{enumerate}

\end{lemma}
\noindent {\bf Proof.} First, consider (a). Since $(x^\star,y^\star)$ is an optimal solution of (P), it satisfies $x^\star\in X$, $y^\star \in Y$, $H(x^\star,y^\star) \in D$, and 
\[
g(x^\star,y^\star) \leq \inf_{z\in Y} \big\{g(x^\star,z)~\big|~H(x^\star,z) \in D\big\} + \tau;
\]
cf. \eqref{eqn:level-reform}. By calmness and Proposition \ref{prop:char_calm}, there is $\lambda^\star \in [0,\infty)$ such that 
\[
\inf_{z\in Y} \big\{g(x^\star,z)~\big|~H(x^\star,z) \in D\big\} = \mu(x^\star,\lambda^\star). 
\]
Consequently, $g(x^\star,y^\star) \leq \mu(x^\star,\lambda^\star) + \tau$ and 
\begin{equation}\label{eqn:withCalmBnd}
\inf \phi_\tau \leq \phi_\tau(x^\star,y^\star,0,0,\lambda^\star) = f(x^\star,y^\star) = \mathfrak{m} \in\reals. 
\end{equation}
We observe that $\inf \phi_\tau \geq \mathfrak{m}$ because of \eqref{eqn:level-reform} and the fact that $\mu(x,\lambda) \leq \inf_{z\in Y} \{g(x,y)$ $|~H(x,y) \in D\}$ for all $x\in\reals^n$ and $\lambda \in [0,\infty)$. Thus, in view of \eqref{eqn:withCalmBnd}, we conclude that $\inf \phi_\tau = \mathfrak{m} \in \reals$.

Second, consider (b). Since $(\bar x, \bar y, \bar u, \bar \alpha, \bar \lambda) \in \nargmin \phi_\tau$, one has $\phi_\tau(\bar x, \bar y, \bar u, \bar\alpha, \bar\lambda) <\infty$. Thus, $\bar x \in X$, $\bar y \in Y$, $\bar u = 0$, $\bar \alpha = 0$, $H(\bar x, \bar y) \in D$, and 
\[
g(\bar x, \bar y) \leq \mu(\bar x, \bar \lambda) + \tau \leq \inf_{z\in Y} \big\{g(\bar x, z)~\big|~H(\bar x, z) \in D\big\} + \tau.
\] 
In view of \eqref{eqn:level-reform}, this means that $(\bar x, \bar y)$ is feasible in (P).  
Moreover, 
\[
\phi_\tau(\bar x, \bar y, \bar u, \bar\alpha, \bar \lambda) = f(\bar x, \bar y) \leq \phi_\tau(x,y,u,\alpha,\lambda) ~~ \forall (x,y,u,\alpha,\lambda)\in \reals^{n+m+q+2}. 
\] 
In particular, consider $x=x^\star$, $y = y^\star$, $u = 0$, $\alpha = 0$, and $\lambda = \lambda^\star\in [0,\infty)$, a value that satisfies
\[
\inf_{z\in Y} \big\{g(x^\star,z)~\big|~H(x^\star,z) \in D\big\} = \mu(x^\star,\lambda^\star).
\]
Such $\lambda^\star$ exists by Proposition \ref{prop:char_calm} and the calmness of the value function $V$ at $x^\star$. This choice yields $f(\bar x, \bar y) \leq \phi_\tau(x^\star ,y^\star , 0, 0,\lambda^\star) = f(x^\star,y^\star) = \mathfrak{m}$ because $x^\star\in X$, $y^\star \in Y$, $H(x^\star,y^\star) \in D$, and 
\[
g(x^\star,y^\star) \leq \inf_{z\in Y} \big\{g(x^\star,z)~\big|~H(x^\star,z) \in D\big\} + \tau = \mu(x^\star,\lambda^\star) + \tau. 
\]
Thus, $(\bar x, \bar y)$ is optimal for (P).

Third, consider (c). Let $\lambda \in [\lambda^\star, \infty)$. By calmness and Proposition \ref{prop:char_calm}, one has
\[
\inf_{z\in Y} \big\{g(x^\star,z)~\big|~H(x^\star,z) \in D\big\} = \mu(x^\star,\lambda). 
\]
Consequently, $g(x^\star,y^\star) \leq \inf_{z\in Y} \{g(x^\star,z)~|~H(x^\star,z) \in D\} + \tau = \mu(x^\star,\lambda) + \tau$ and $\phi_\tau(x^\star, y^\star, 0, 0, \lambda)$ $=$ $f(x^\star,y^\star) = \mathfrak{m} = \inf \phi_\tau$ by (a).\eop

\noindent {\bf Proof of Theorem \ref{thm:withcalmness}(a).} We invoke Lemma \ref{lemma:limsup}, with $\beta = \tau$, and mimic the beginning part of the proof of 
Theorem \ref{thm:withoutcalmness}(a). This leads to $\nlimsup \mathfrak{m}^\nu = \nlimsup(\inf \phi^\nu) \leq \inf \phi_{\tau} = \mathfrak{m} \in \reals$, 
where the last equality follows from Lemma \ref{lemma:phi}(a).\eop

\noindent {\bf Proof of Theorem \ref{thm:withcalmness}(b).} Lemmas \ref{lemma:liminf} and \ref{lemma:limsup} hold with $\beta = \tau$. This implies that $\phi^\nu \eto \phi_\tau$; see \eqref{eqn:liminf} and \eqref{eqn:limsup}. The existence of a compact set $B$ implies that $\{\phi^\nu, \nu\in\nats\}$ is tight. We bring in \cite[Theorem 5.5(d)]{primer} to conclude that $\inf \phi^\nu \to \inf \phi_\tau$ whenever $\inf \phi_\tau<\infty$ and $\phi_\tau$ takes nowhere the value $-\infty$. The former condition holds by Lemma \ref{lemma:phi}(a). The latter condition holds because $f(x,y)>-\infty$ for all $x\in\reals^n$ and $y\in \reals^m$. This confirms that $\mathfrak{m}^\nu \to \mathfrak{m} \in \reals$ after invoking Lemma \ref{lemma:phi}(a). 

We observe that $(x^\star, y^\star, 0, 0, \lambda) \in \nargmin \phi_\tau$ by Lemma \ref{lemma:phi}(c). Moreover, \cite[Theorem 5.5(e)]{primer} guarantees the existence of vanishing $\{\epsilon^\nu \in [0,\infty), \nu\in\nats\}$ and $\{w^\nu\in \reals^{n+m + q +2}, \nu\in\nats\}$ such that $w^\nu \in \epsilon^\nu$-$\nargmin \phi^\nu$ and $w^\nu \to (x^\star, y^\star, 0, 0, \lambda)$. Since the $\epsilon^\nu$-optimal solutions of (P)$^\nu$ coincides with $\epsilon^\nu$-$\nargmin \phi^\nu$, the conclusion holds.\eop

\noindent {\bf Proof of Theorem \ref{thm:withcalmness}(c).} Lemmas \ref{lemma:liminf} and \ref{lemma:limsup} hold with $\beta = \tau$. This implies that $\phi^\nu \eto \phi_\tau$. Moreover, $\inf \phi_\tau<\infty$ holds by Lemma \ref{lemma:phi}(a) and $\phi_\tau$ takes nowhere the value $-\infty$ trivially. Thus, \cite[Theorem 5.5(b)]{primer} applies and yields $(\hat x, \hat y, \hat u, \hat \alpha, \hat \lambda) \in \nargmin \phi_\tau$. 
Lemma \ref{lemma:phi}(b) ensures that $(\hat x, \hat y)$ is an optimal solution of (P). By \cite[Theorem 5.5(c)]{primer}, $\inf \phi^\nu \Nto \inf \phi_\tau$. By Lemma \ref{lemma:phi}(a), $\inf \phi_\tau = \mathfrak{m}\in\reals$ and the conclusion follows.\eop

\noindent {\bf Proof of Proposition \ref{prop:suffCalm}.} If $Y$ or $D$ is empty, the conclusion follows trivially. Thus, we concentrate on the case with nonempty $Y$ and $D$. We note that $S(x)$ is nonempty because otherwise \eqref{eqn:regularity} would not hold. Let $\epsilon \in (0,\infty)$ and $y'\in Y$. Then, there exists $y \in S(x)$ such that $\|y - y'\|_2 \leq \dist_2(y', S(x)) + \epsilon$. By \eqref{eqn:LipAssump} and \eqref{eqn:regularity}, we obtain 
\begin{align*}
g(x,y) & \leq g(x,y') + \kappa_1\|y - y'\|_2 \leq g(x,y') + \kappa_1 \dist_2\big(y', S(x)\big) + \kappa_1 \epsilon\\
& \leq g(x,y') + \kappa_1 \kappa_2 \dist\big(H(x,y'),D\big) + \kappa_1 \epsilon.  
\end{align*}
Since $y'\in Y$ is arbitrary, one has $g(x,y) \leq \inf_{z\in Y} g(x,z) + \kappa_1 \kappa_2 \dist(H(x,z),D) + \kappa_1 \epsilon$. Moreover, $y$ satisfies $H(x,y)\in D$ so that $\inf_{z\in Y}\{g(x,z) ~|~H(x,z) \in D\} \leq g(x,y)$. Consequently,  
\[
V(x,0) \leq \inf_{z\in Y} g(x,z) + \kappa_1 \kappa_2 \dist\big(H(x,z),D\big) + \kappa_1 \epsilon. 
\]
The last term drops off because $\epsilon$ is arbitrary. Also, $V(x,0) \geq \inf_{z\in Y} g(x,z) + \kappa_1 \kappa_2 \dist(H(x,z),D)$ and the conclusion follows via Proposition \ref{prop:char_calm}.\eop

\noindent {\bf Proof of Theorem \ref{thm:withLocalCalmness}.} Lemmas \ref{lemma:liminf} and \ref{lemma:limsup} hold with $\beta = \tau$. This implies that $\phi^\nu \eto \phi_\tau$. We observe that $\phi_\tau$ nowhere takes the value $-\infty$ and $\inf \phi_\tau<\infty$ by Lemma \ref{lemma:phi}(a). Thus, \cite[Theorem 5.5(b)]{primer} applies and yields that every cluster point of a sequence of $\omega^\nu$-minimizers of $\phi^\nu$ is contained in $\nargmin \phi_\tau$ provided that $\omega^\nu\to 0$. Along any subsequence of such convergent $\omega^\nu$-minimizers, \cite[Theorem 5.5(c)]{primer} ensures that $\inf \phi^\nu$ converges to $\inf \phi_\tau$.  

By local calmness of $V$ at $\hat x$, there exist $\lambda\in [0,\infty)$ and $\rho \in (0,\infty)$ such that 
  \[
  \mu(x,\lambda) = \inf_{z\in Y} g(x,z) + \lambda \dist\big(H(x,z), D\big) = \inf_{z\in Y} \big\{g(x,z) ~\big|~ H(x,z) \in D\big\} ~~~~\forall x\in \ball(\hat x, \rho) 
  \]
as seen by Proposition \ref{prop:char_calm} and Definition \ref{def:localcalm}. Since $\bar \lambda^\nu\to \infty$ by Assumption \ref{ass:basic}(f), and $x^\nu\Nto \hat x$, there exists $\bar \nu$ such that $x^\nu \in \ball(\hat x, \rho)$ and $\bar\lambda^\nu \geq \lambda$ for all $\nu\geq \bar \nu$, $\nu\in N$. 

For $\nu\in N$, we construct 
\begin{align*}
\tilde \alpha^\nu & = \alpha^\nu - \mu^\nu\big(x^\nu, \max\{\lambda, \lambda^\nu\}\big) + \mu^\nu(x^\nu,\lambda)\\
\omega^\nu & = \epsilon^\nu + \theta^\nu \Big(\mu^\nu\big(x^\nu, \max\{\lambda, \lambda^\nu\}\big) - \mu^\nu(x^\nu,\lambda)  \Big).
\end{align*}
We claim that $(x^\nu, y^\nu, u^\nu, \tilde \alpha^\nu, \lambda)$ is an $\omega^\nu$-optimal solution for (P)$^\nu$ for all $\nu\geq \bar \nu$, $\nu\in N$, and $\omega^\nu\Nto 0$ and $\tilde \alpha^\nu \Nto \hat \alpha$. To confirm this claim, we proceed in three steps:

\smallskip

\noindent {\bf Step 1.} Neither $\mu^\nu(x^\nu, \max\{\lambda, \lambda^\nu\})$ nor $\mu^\nu(x^\nu,\lambda)$ is equal to $\infty$ because $Y^\nu$ and $D$ are nonempty. Since $(x^\nu, y^\nu, u^\nu, \alpha^\nu, \lambda^\nu)$ is an $\epsilon^\nu$-optimal solution of (P)$^\nu$, it must be feasible. In turn, this implies that \eqref{eqn:PnuConstr_mu} holds at $(x^\nu, y^\nu, u^\nu, \alpha^\nu, \lambda^\nu)$  and thus $- \infty < \mu^\nu(x^\nu,\lambda^\nu) \leq \mu^\nu(x^\nu, \max\{\lambda, \lambda^\nu\})$, where the last inequality follows by the nondecreasing property of $\mu^\nu(x^\nu,\cdot)$. 
We also confirm that $\mu^\nu(x^\nu,\lambda)>-\infty$ for $\nu\geq \bar \nu$, $\nu\in N$, and in the process achieve a bound on the distance between $\mu^\nu(x^\nu,\lambda)$ and $\mu^\nu( x^\nu, \max\{\lambda, \lambda^\nu\})$. We consider settings (a) and (b) in the theorem separately. 

Under (a): Without loss of generality, we assume that $Y^\nu = Y$ for all $\nu \geq \bar\nu$. For $\nu\geq \bar \nu$, $\nu\in N$, 
\begin{align*}
&\mu^\nu\big(x^\nu, \max\{\lambda, \lambda^\nu\}\big) = \inf_{z\in Y} g^\nu(x^\nu,z) + \max\{\lambda, \lambda^\nu\} \dist\big(H^\nu(x^\nu, z), D\big)\\     
&\leq \inf_{z\in Y} g(x^\nu,z) + \max\{\lambda, \lambda^\nu\} \dist\big(H(x^\nu, z), D\big) + \delta^\nu + \bar\lambda^\nu \eta^\nu\\     
& = \mu(x^\nu,\lambda) + \delta^\nu + \bar\lambda^\nu \eta^\nu\\     
& \leq \inf_{z\in Y} g^\nu(x^\nu,z) + \lambda \dist\big(H^\nu(x^\nu, z), D\big) + 2(\delta^\nu + \bar\lambda^\nu \eta^\nu) = \mu^\nu(x^\nu,\lambda) + 2(\delta^\nu + \bar\lambda^\nu \eta^\nu), 
\end{align*}
where we use Assumption \ref{ass:basic}(c,d) and the fact that $\mu(x^\nu, \max\{\lambda, \lambda^\nu\}) = \mu(x^\nu,\lambda)$ because $\lambda$ is a penalty threshold for $V$ at $x^\nu$. Since $\mu^\nu(x^\nu, \max\{\lambda, \lambda^\nu\}) > - \infty$, this yields $\mu^\nu( x^\nu,\lambda)>-\infty$ for $\nu\geq \bar \nu$, $\nu\in N$.

Under (b): Suppose that $\nu\geq \bar \nu$, $\nu\in N$. Without loss of generality we assume that $\sup_{y\in Y} \dist_2(y,Y^\nu)$ $< \infty$ because this term vanishes as $\nu\to\infty$. Let $\gamma \in (0,\infty)$. We consider two cases. First, suppose that $\mu(x^\nu, \max\{\lambda, \lambda^\nu\}) = -\infty$. Then, there exists $\bar y \in Y$ such that 
\[
g(x^\nu,\bar y) + \max\{\lambda,  \lambda^\nu\} \dist\big(H( x^\nu, \bar y), D\big) \leq -1/\gamma
\]
and there is $\bar y^\nu \in Y^\nu$ such that $\|\bar y^\nu - \bar y\|_2 \leq \dist_2(\bar y,Y^\nu) + \gamma$. Using Assumption \ref{ass:basic}(c,d), we obtain
\begin{align*}
&\mu^\nu\big(x^\nu, \max\{\lambda, \lambda^\nu\}\big)\\
&\leq \inf_{z\in Y^\nu} g(x^\nu,z) + \max\{\lambda, \lambda^\nu\} \dist\big(H(x^\nu, z), D\big) + \delta^\nu + \bar\lambda^\nu \eta^\nu\\     
&\leq g(x^\nu,\bar y^\nu) + \max\{\lambda, \lambda^\nu\} \dist\big(H(x^\nu, \bar y^\nu), D\big) + \delta^\nu + \bar\lambda^\nu \eta^\nu\\     
&\leq g(x^\nu,\bar y) + \max\{\lambda,  \lambda^\nu\} \dist\big(H( x^\nu, \bar y), D\big) + (\kappa_g + \bar\lambda^\nu\kappa_H)\|\bar y^\nu-\bar y\|_2 + \delta^\nu + \bar\lambda^\nu \eta^\nu\\     
&\leq -1/\gamma + (\kappa_g + \bar\lambda^\nu\kappa_H)\|\bar y^\nu-\bar y\|_2 + \delta^\nu + \bar\lambda^\nu \eta^\nu\\ 
&\leq -1/\gamma + (\kappa_g + \bar\lambda^\nu\kappa_H)\sup_{y\in Y} \dist_2(y,Y^\nu) + \delta^\nu + \bar\lambda^\nu \eta^\nu + \gamma.
\end{align*}     
The terms after $-1/\gamma$ are finite. Since $\gamma$ is arbitrary, this implies that the right-hand side is an arbitrarily low number. Thus, we have reached a contradiction because $\mu^\nu(x^\nu, \max\{\lambda, \lambda^\nu\})>-\infty$. 

Second, suppose that $\mu(x^\nu, \max\{\lambda, \lambda^\nu\})$ is finite. Then, there is $\bar y \in Y$ such that 
\begin{align*}
&g(x^\nu,\bar y) + \max\{\lambda, \lambda^\nu\} \dist\big(H(x^\nu, \bar y), D\big) \\
&\leq \inf_{z\in Y} g(x^\nu,z) + \max\{\lambda, \lambda^\nu\} \dist\big(H(x^\nu, z), D\big) + \gamma
\end{align*}
and  there is $\bar y^\nu \in Y^\nu$ such that $\|\bar y^\nu - \bar y\|_2 \leq \dist_2(\bar y,Y^\nu) + \gamma$. Again, Assumption \ref{ass:basic}(c,d) yields
\begin{align*}
&\mu^\nu\big(x^\nu, \max\{\lambda, \lambda^\nu\}\big) \leq \inf_{z\in Y^\nu} g(x^\nu,z) + \max\{\lambda, \lambda^\nu\} \dist\big(H(x^\nu, z), D\big) + \delta^\nu + \bar\lambda^\nu \eta^\nu\\     
&\leq g(x^\nu,\bar y^\nu) + \max\{\lambda, \lambda^\nu\} \dist\big(H(x^\nu, \bar y^\nu), D\big) + \delta^\nu + \bar\lambda^\nu \eta^\nu\\     
&\leq g(x^\nu,\bar y) + \max\{\lambda, \lambda^\nu\} \dist\big(H(x^\nu, \bar y), D\big) + (\kappa_g + \bar\lambda^\nu\kappa_H)\|\bar y^\nu-\bar y\|_2 + \delta^\nu + \bar\lambda^\nu \eta^\nu\\     
&\leq \inf_{z\in Y} g(x^\nu,z) + \max\{\lambda, \lambda^\nu\} \dist\big(H( x^\nu, z), D\big) + (\kappa_g + \bar\lambda^\nu\kappa_H)\|\bar y^\nu-\bar y\|_2 + \delta^\nu\\
&~~~ + \bar\lambda^\nu \eta^\nu + \gamma\\     
&\leq \mu(x^\nu,\lambda) + (\kappa_g + \bar\lambda^\nu\kappa_H)\sup_{y \in Y} \dist_2(y,Y^\nu)  + \delta^\nu + \bar\lambda^\nu \eta^\nu + 2\gamma\\     
&\leq \inf_{z\in Y} g^\nu( x^\nu,z) + \lambda \dist\big(H^\nu(x^\nu, z), D\big) + (\kappa_g + \bar\lambda^\nu\kappa_H)\sup_{y \in Y} \dist_2(y,Y^\nu)\\
&~~~  + 2(\delta^\nu + \bar\lambda^\nu \eta^\nu) + 2\gamma\\     
&\leq  \mu^\nu(x^\nu,\lambda) + (\kappa_g + \bar\lambda^\nu\kappa_H)\sup_{y \in Y} \dist_2(y,Y^\nu)  + 2(\delta^\nu + \bar\lambda^\nu \eta^\nu) + 2\gamma.  
\end{align*}
Since $\mu^\nu(x^\nu, \max\{\lambda, \lambda^\nu\}) > - \infty$, this yields $\mu^\nu(x^\nu,\lambda)>-\infty$ for $\nu\geq \bar \nu$, $\nu\in N$ because the other terms on the right-hand side are finite. In fact, since $\gamma$ is arbitrary, 
\begin{align}
&0 \leq \mu^\nu\big(x^\nu, \max\{\lambda, \lambda^\nu\}\big) - \mu^\nu( x^\nu,\lambda)\label{eqn:28munu}\\
& \leq (\kappa_g + \bar\lambda^\nu\kappa_H)\sup_{y \in Y} \dist_2(y,Y^\nu)  + 2(\delta^\nu + \bar\lambda^\nu \eta^\nu).\nonumber  
\end{align}

\noindent {\bf Step 2.} Since $\mu^\nu(x^\nu, \max\{\lambda, \lambda^\nu\})$ and $\mu^\nu( x^\nu,\lambda)$ are finite by Step 1, $\omega^\nu$ and $\tilde \alpha^\nu$ are finite. In setting (a), also by Step 1, we have $0 \leq \mu^\nu\big(x^\nu, \max\{\lambda, \lambda^\nu\}\big) - \mu^\nu(x^\nu,\lambda) \leq 2(\delta^\nu + \bar\lambda^\nu \eta^\nu)$ for all $\nu\geq \bar\nu$, $\nu\in N$. Invoking Assumption \ref{ass:basic}(f), this implies 
\[
\theta^\nu \Big(\mu^\nu\big(x^\nu, \max\{\lambda, \lambda^\nu\}\big) - \mu^\nu(x^\nu,\lambda)  \Big) \leq 2 \theta^\nu(\delta^\nu + \bar\lambda^\nu \eta^\nu) \Nto 0. 
\]
In setting (b), we similarly obtain via \eqref{eqn:28munu} and Assumption \ref{ass:basic}(f) that 
\begin{align*}
&\theta^\nu \Big(\mu^\nu\big(x^\nu, \max\{\lambda, \lambda^\nu\}\big) - \mu^\nu(x^\nu,\lambda)  \Big)\\
& \leq \theta^\nu\big((\kappa_g + \bar\lambda^\nu\kappa_H)\sup_{y \in Y} \dist_2(y,Y^\nu)  + 2(\delta^\nu + \bar\lambda^\nu \eta^\nu)\big) \Nto 0. 
\end{align*}
In either setting, we conclude that $\omega^\nu \Nto 0$ because $\epsilon^\nu\to 0$. Moreover, $\tilde \alpha^\nu \Nto \hat \alpha$ because $\theta^\nu\to \infty$.

\smallskip

\noindent {\bf Step 3.} The objective function value in (P)$^\nu$ at the point $(x^\nu, y^\nu, u^\nu, \tilde \alpha^\nu, \lambda)$ is 
\begin{align*}
& f^\nu(x^\nu, y^\nu) + \sigma^\nu\|u^\nu\| - \theta^\nu \tilde \alpha^\nu\\ 
& = f^\nu(x^\nu, y^\nu) + \sigma^\nu\|u^\nu\| - \theta^\nu \alpha^\nu + \theta^\nu \Big(\mu^\nu\big(x^\nu, \max\{\lambda, \lambda^\nu\}\big) - \mu^\nu( x^\nu,\lambda)  \Big)\\
& \leq \mathfrak{m}^\nu + \epsilon^\nu + \theta^\nu \Big( \mu^\nu\big(x^\nu, \max\{\lambda, \lambda^\nu\}\big) - \mu^\nu( x^\nu,\lambda)  \Big) = \mathfrak{m}^\nu + \omega^\nu.
\end{align*}
We also check the constraints \eqref{eqn:PnuConstr_mu}: 
\begin{align*}
g^\nu( x^\nu,  y^\nu) + \tilde \alpha^\nu &\leq g^\nu( x^\nu,  y^\nu) +  \alpha^\nu - \mu^\nu\big( x^\nu, \max\{\lambda,  \lambda^\nu\}\big) + \mu^\nu( x^\nu,\lambda)\\ 
&\leq \mu^\nu( x^\nu,  \lambda^\nu) + \tau^\nu - \mu^\nu\big( x^\nu, \max\{\lambda,  \lambda^\nu\}\big) + \mu^\nu( x^\nu,\lambda)\\ 
&\leq \mu^\nu\big( x^\nu, \max\{\lambda,  \lambda^\nu\}\big) + \tau^\nu - \mu^\nu\big( x^\nu, \max\{\lambda,  \lambda^\nu\}\big) + \mu^\nu( x^\nu,\lambda)\\ 
&\leq \mu^\nu( x^\nu,\lambda) + \tau^\nu. 
\end{align*}
Thus, \eqref{eqn:PnuConstr_mu} holds for $( x^\nu,  y^\nu,  u^\nu, \tilde \alpha^\nu, \lambda)$. Moreover, $\lambda \in [0, \bar\lambda^\nu]$ and $\tilde \alpha^\nu \leq  \alpha^\nu\leq 0$ for $\nu\geq \bar \nu, \nu\in N$. We have confirmed that $( x^\nu, y^\nu,  u^\nu, \tilde \alpha^\nu, \lambda)$ is an $\omega^\nu$-optimal solution for (P)$^\nu$.

Under both (a) and (b), $(\hat x, \hat y, \hat u, \hat \alpha, \lambda) \in \nargmin \phi_\tau$. It follows by Lemma \ref{lemma:phi}(b) that $(\hat x, \hat y)$ is an optimal solution for (P). Since $\inf \phi^\nu \Nto \phi_\tau$ by \cite[Theorem 5.5(c)]{primer}, Lemma \ref{lemma:phi}(a) yields $\mathfrak{m}^\nu \Nto \mathfrak{m}$.\eop

\vspace{0.2cm}

\state Acknowledgement. This work is supported in part by the Office of Naval Research under grants N00014-24-1-2277, N00014-24-1-2318, and N00014-24-1-2492. 

\bibliographystyle{plain}
\bibliography{refs}

\end{document}